\documentclass[11pt,reqno]{amsart}

\newtheorem{theorem}{Theorem}[section]
\newtheorem{lemma}[theorem]{Lemma}
\newtheorem{corollary}[theorem]{Corollary}

\theoremstyle{definition}
\newtheorem{assumption}[theorem]{Assumption}
\newtheorem{definition}[theorem]{Definition}
  \newtheorem{exercise*}[theorem]{Exercise*}

\theoremstyle{remark} 
\newtheorem{example}[theorem]{Example}
\newtheorem{remark}[theorem]{Remark}
 
\newcommand{\loc}{\textnormal{loc}}
\newcommand{\shharp}{=\kern -.5em\|}
\newcommand{\vsharp}{\asymp\kern -.5em\|}

\newcommand{\infsup}{\operatornamewithlimits{\,inf\,\,\,\,sup}}

\newcommand\bM{\mathbb{M}}
\newcommand\bR{\mathbb{R}}

\newcommand\bS{\mathbb{S}}

\newcommand\chB{\check B}
\newcommand\chC{\check C}

\newcommand\sfu{{\sf u}} 
\newcommand\sfv{{\sf v}}

\newcommand\cL{\mathcal{L}}

\newcommand{\sign}{\text{\rm\,sign}\,}

\newcommand{\osc}{\operatornamewithlimits{osc}}

\newcommand{\diam}{{\rm diam}}

\makeatletter
\def\dashintindex{\operatorname%
{-\kern-.7em\DOTSI\intop\ilimits@}}%
\def\dashint{\operatorname%
{\,\,\text{\bf--}\kern-.98em\DOTSI\intop\ilimits@\!\!}}
\makeatother

 \newcommand{\mysection}[1]{\section{#1}
 \setcounter{equation}{0}}

\renewcommand{\eqref}[1]{\text{\rm(\ref{#1})}}

\renewcommand\({{\rm(}}

\renewcommand\){{\rm)}}

\begin{document}  

\title[Fully nonlinear
parabolic 
equations]
{On the existence of $W^{1,2}_{p}$ 
solutions for fully nonlinear parabolic
equations under either relaxed 
or no  convexity assumptions}

\author{N.V. Krylov}

\email{nkrylov@umn.edu}
\address{127 Vincent Hall, University of Minnesota,
 Minneapolis, MN, 55455}

 \keywords{Fully nonlinear
parabolic  
equations, cut-off equations}

\subjclass{35K55, 35K20}

\begin{abstract}
We establish the existence  of solutions
of fully nonlinear parabolic second-order equations
like $\partial_{t}u+H(v,Dv,D^{2}v,t,x)=0$
in smooth cylinders without requiring $H$
to be convex or concave with respect to the second-order derivatives.
Apart from ellipticity nothing is required of $H$ at points at which
 $|D^{2}v|\leq K$,
where $K$ is any fixed constant. For large $|D^{2}v|$ some kind
of relaxed convexity assumption with respect to $D^{2}v$ mixed 
with a VMO condition
with respect to $t,x$ are still imposed. The solutions are sought
in Sobolev classes. We also establish the solvability
without almost any conditions on $H$, apart from ellipticity,
but of a ``cut-off'' version of the equation 
$\partial_{t}u+H(v,Dv,D^{2}v, t, x)=0$.
\end{abstract}

\maketitle

\mysection{Introduction and main results}
                                        \label{section 6.21.02}

In this paper
 we consider parabolic  equations
\begin{equation}
                                                \label{8,17.7}
\partial_{t}v(t,x)+H[v](t,x)  =0,
\end{equation}
where 
$$
H[v](t,x)= H(v(t, x),D v(t, x),D^{2}v(t, x),t ,x),
$$
in subdomains of 
$$
\bR^{d +1}=\{(t,x):t\in\bR,x\in\bR^{d}\}.
$$ 

Let $\Omega\in C^{1,1}$ be an open bounded subset of $\bR^{d}$.
  Fix $T\in(0,\infty)$
and set
$$
\Pi=[0,T)\times\Omega
$$
(if the $t$-axis is directed vertically, $[0,T)\times\Omega$
is indeed looking like a pie).
Fix
  $$
p>d   \quad\text{and a measurable function}\quad
  \bar G\geq0 \quad\text{on}\quad \bR^{d+1}.
$$

One of our
 main results implies that, for $d=3$,    equation
\eqref{8,17.7} with
($a\wedge b=\min(a,b)$)
$$
H(D^{2}u,x):=\bar G(t,x)\wedge|D_{12}u|
+\bar G(t,x)\wedge|D_{23}u|+\bar G(t,x)\wedge|D_{31}u|
$$
\begin{equation}
                                                      \label{3.4.1}
+ \Delta u-f(t,x) 
\end{equation}
in $\Pi$ with zero boundary condition on its parabolic boundary has a unique
solution $u\in W^{1,2}_{p}(\Pi)$, provided that $\bar G,f\in L_{p}(\Pi)$
with $p>d+2$.
Recall that $W^{1,2}_p(\Pi)$ denotes the set of
functions $v$ defined in
$\Pi$ such that $\partial_{t}v$, $v$, $Dv
 $, and $D^2v $ 
are in
$L_p(\Pi)$. Observe that $H$ in \eqref{3.4.1}
is neither convex nor concave with respect to $D^{2}u$.
So far, there are only two approaches to such equations:
the theory of ($L_{p}$) viscosity solutions
and the theory of stochastic differential games,
provided $H$ has a somewhat special form. The past experience
shows that
it is hard
to expect getting sharp quantitative results using probability theory.
On the other hand, the theory of viscosity solutions indeed produced
some remarkable quantitative results (see, for instance,
\cite{CKS00}, \cite{JS_05}  and the references therein).
However, to the best of the author's knowledge the result
stated above about \eqref{3.4.1}
is either very hard to obtain by using the theory
of ($L_{p}$) viscosity solutions or is just beyond it,
at least at the current stage.
It seems that the best information, that theory provides
at the moment,
is the existence of the maximal and minimal
$(L_{p})$ viscosity solution (see \cite{JS_05}),
no uniqueness of $(L_{p})$ viscosity solutions
can be inferred  for \eqref{3.4.1} and no regularity
apart from the classical $C^{\alpha}$-regularity (see
\cite{CKLS99}). 

The current paper is a natural continuation
of \cite{Kr_13}
where similar results are obtained for elliptic equations.

    Fix some constants $K_{0},K_{F}\in[0,\infty)$,
$\delta\in(0,1]$.
Denote by $\bS$ the set of symmetric $d\times d$ matrices
and let $\bS_{\delta}$ be the subset of $\bS$
consisting of matrices $a$ such that
$$
\delta|\lambda|^{2}\leq a^{ij}\lambda^{i}\lambda^{j}
\leq \delta^{-1}|\lambda|^{2}\quad \forall \lambda\in\bR^{d}.
$$

Here are our assumptions about $H$.

  \begin{assumption}
                                    \label{assumption 7,16.1} 

The function   $H(\sfu ,t, x)$,
$$
\sfu=(\sfu',\sfu''),\quad \sfu'=(\sfu'_{0},\sfu'_{1},...,\sfu'_{d})
\in \bR^{d+1},\quad
\sfu''\in\bS,\quad (t,x)\in\bR^{d+1},
$$
is measurable with respect to $(\sfu', t,x)$.

 \end{assumption}

The following assumptions contain (small)  parameters $\hat\theta,\theta\in(0,1]$
which are specified later in our results.

\begin{assumption}
                                    \label{assumption 7,16.01}
There are two measurable
functions 
$$
F(\sfu,t,x)=F(\sfu'_{0},\sfu'',t,x),\quad G(\sfu,t,x)
$$
 such that
$$
H =F +G .
$$
 For $\sfu''\in\bS,\sfu'\in\bR^{d+1}$, and $(t,x)\in\bR^{d+1}$ we have
$$
|G( \sfu ,t,x)|\leq \hat\theta|\sfu''|+ K_{0}|\sfu'|+\bar{G}(t,x),\quad 
F(0,t,x)\equiv0.
$$
\end{assumption}

Define
$$
B_{R}(x_{0})=\{x\in\bR:|x-x_{0}|<R\},\quad B_{R}=B_{R}(0),
$$
$$
C_{r}(t_{0},x_{0})=[t_{0},t_{0}+r^{2})\times B_{r}(x_{0}),
\quad C_{r}=C_{r}(0),
$$
and for Borel $\Gamma\subset \bR^{d+1}$ denote
\index{$|\Gamma$}%
 by $|\Gamma|$
the volume of $\Gamma$.

\begin{assumption}
                                \label{assumption 8,18.2} 
  (i) The function $F$ is Lipschitz continuous with respect to $\sfu''$
with Lipschitz constant $K_{F}$ and is measurable
with respect to $(t,x)$.

Moreover there exist  $R_0\in(0,1]$ and $\tau_{0}\in[0,\infty)$
such that, for any $\sfu'_{0}\in\bR$,
 $z_{0}=(t_{0},x_{0})\in\Pi$ and $r\in (0, R_0]$,
one can find a {\em convex\/} function $\bar{F} (\sfu'' )=
\bar{F}_{z_{0},r,\sfu'_{0}} (\sfu'' )$ (depending only on $\sfu''$)  for which

(ii\,) We have $\bar{F}(0)=0$ and at all points of differentiability
of $\bar{F}$ we have $D_{\sfu''}\bar{F} \in\bS_{\delta}$;
 
(iii\,) 
For any $\sfu''\in\bS$ with $|\sfu''|=1$,  
 we have
\begin{equation}
                                                \label{8,18.3}
\int_{ \hat C_{r}(z_{0})}\sup_{\tau>\tau_{0}}\tau^{-1}
|F(\sfu'_{0},\tau \sfu'' ,z)-\bar{F}(\tau \sfu'')| \,dz \leq \theta
|\hat C_{r}(z_{0})|,
\end{equation}
where
$$
\hat C_{r}(z_{0})=(t_{0},t_{0}+r^{2})
\times (\Omega\cap B_{r}(x_{0}));
$$

 (iv) There exists a continuous increasing 
function $\omega_{F }(\tau)$, $\tau\geq0$,
such that $\omega_{F}(0)=0$ 
and for any $\sfu'_{0},\sfv'_{0}\in\bR$, $(t,x)
\in\Pi$, and $\sfu''\in\bS$ we have
$$
|F(\sfu'_{0},  \sfu'' ,t,x)-F(\sfv'_{0},  \sfu'',t ,x)|\leq
 \omega_{F}(|\sfu'_{0}-\sfv'_{0}|)|u''|.
$$ 

\end{assumption}

\begin{remark}
                                            \label{remark 1,19,2}
Assumptions \ref{assumption 7,16.01} and
 \ref{assumption 8,18.2} (iv)  imply that
\begin{equation}
                                                   \label{1,19,500}
 |H(\sfu',0,t,x)| \leq K_{0}|\sfu'|+\bar{G}(t,x)\quad\forall \sfu',(t,x)
\in\bR^{d+1}.
\end{equation}

\end{remark}

\begin{assumption}
                                \label{assumption 8,18.3} 
We are given a function $g\in W^{1,2}_{p}(\Pi)$.
\end{assumption}

If $z_{i}=(t_{i},x_{i})\in\bR^{d+1}$, $i=1,2$, we set
$$
\rho(z_{1},z_{2})=|t_{1}-t_{2}|^{1/2}+|x_{1}-x_{2}|.
$$

\begin{definition}
                                          \label{definition 1,18,1}
For a function $u\in C(\bar\Pi)$ 
set
$$
\omega_{u}(\Pi,\rho)=
\sup \{
 |u(z_{1})-u(z_{2})|  :z_{1},z_{2}\in\Pi,
\rho(z_{1},z_{2}) \leq\rho\},
$$
$$
\omega_{F,u,\Pi}(\rho)=\omega_{F}(\omega_{u}(\Pi,\rho)).
$$

  We will sometimes say that a certain constant 
depends only on A,B,..., and
 {\em the function\/}
$\omega_{F,u,\Pi}$. This is to mean
that it depends only on A,B,..., and on the maximal
 solution of an inequality like  $N_{0}\omega_{F,u,\Pi}(\rho)\leq
1/2$, where the range of $\rho$ and
 the value of $N_{0}$ depending only on 
A,B,... could be always found out from our arguments.

\end{definition}
In the following theorem about a priori estimates
there is no ellipticity assumption on $H$.
If $Q$ is a subdomain in $\bR^{d+1}$, by $\partial'Q$
we denote its parabolic boundary.

\begin{theorem}
                                    \label{theorem 8,17.3}
Let $p>d+1$. Then
there exist constants $ \theta,\hat\theta\in (0,1]$, depending
only on $d$, $p$, $\delta$, $K_{F}$,  
and $M_{2}( \Omega)$ \($\rho_{0}(\Omega)$ and
 $M_{2}(\Omega)$ are introduced later\),  
such that, if Assumptions \ref{assumption 7,16.01}
and \ref{assumption 8,18.2}  
are satisfied with these
 $\hat \theta$ and $ \theta$, respectively, then for any
 $u\in W^{1, 2}_{d+1}(\Pi ) $
 that satisfies \eqref{8,17.7} in $\Pi $ \(a.e.\)  
and equals $g$ on $\partial'\Pi$
  we have  
\begin{equation}
                                                \label{8,18.1}
\|u\|_{W^{1,2}_{p}(\Pi )}\le N
\| \bar G\|_{L_p(\Pi )}+N \|g\|_{W^{1,2}_{p}(\Pi )}
 +N \tau_{0}+N\sup_{\Pi}|u| ,
\end{equation}
where  
the constants $N$  depend only on
 $K_{0}$, $K_{F}$, $d$, $p$, $\delta$,  $R_{0}$, $\rho_{0}(\Omega)$,
 $M_{2}(\Omega)$, $\diam(\Omega)$, $T$,
 and the functions $\omega_{F,u,\Pi}$ 
 and $\omega_{F,g,\Pi}$. 
 
\end{theorem}

In the literature,   {\em interior\/} $ W^2_p, p>d$, {\em a priori\/}
estimates for a class of
fully nonlinear uniformly {\em elliptic\/} equations in $\bR^{d}$  in
 the framework of viscosity solutions were first obtained by
 Caffarelli in \cite{Caf89}  
(see also
\cite{CC95}).
 Adapting his technique, similar interior 
a priori estimates were proved by Wang
\cite{Wa92} for parabolic equations. In the same paper, a boundary
estimate is stated but without   proof; see Theorem 5.8 there.  By
exploiting a weak reverse H\"older's inequality, the result of
\cite{Caf89} was sharpened by Escauriaza in \cite{Es93}, 
who obtained the
interior $ W^2_p$-estimate for the same equations allowing
$ p>d-\varepsilon$, with a small constant
$\varepsilon>0$ depending only on the ellipticity constant and $d$. 

The above cited works    are quite remarkable 
in one   respect--they do not suppose that $H$ is convex
 or concave in $ D^{2}u$.
But  they only show that 
to prove a priori estimates
it suffices to prove 
the interior $  C^{2}$--estimates for
``harmonic'' functions. However, up to now,
these estimates  are only known under convexity assumptions.

Also obtaining boundary $ W^{2}_{p}$ estimate by using the theory
of viscosity solutions
 turned out to be  extremely challenging and only
in 2009, twenty years after the work of Caffarelli,
 Winter \cite{Wi09}   proved the solvability
 in $ W^{2}_{p}(\Omega)$ of equations  
 with Dirichlet boundary condition in $ \Omega\in C^{1,1}$. 

It is
also worth noting that a solvability theorem in the space
$  W^{1,2}_{p,\text{ loc}}(\Pi)\cap C(\bar \Pi)$ is given  in
M. G. Crandall, M. Kocan, A. \'Swi{\c e}ch \cite{CKS00} for
the boundary-value problem for fully nonlinear parabolic equations. 
The above mentioned existence results of \cite{CKS00} and \cite{Wi09}
are proved under the assumption that $ H$ is convex in $ D^{2}v$
and in all papers mentioned above
 a small oscillation assumption in the integral sense is
imposed on the operators. In the case of linear equations
this  small oscillation assumption is equivalent to requiring
the main coefficients to be uniformly close to uniformly
continuous ones. Our Assumption \ref{assumption 8,18.2} 
is satisfied in this case if 
 the main coefficients  are just in VMO.
The above cited works  are performed in the framework of
viscosity solutions. 

To the best of the author's knowledge the only
paper treating the solvability in the {\em global\/}
Sobolev spaces for {\em parabolic\/} equations is \cite{DKL_12},
where the assumptions are much heavier than here.

To have the solvability we need ellipticity
and more regularity of $H$.  

\begin{assumption}
                                 \label{assumption 7,16.5}
 For any $(t,x)\in\bR^{d+1}$, the function   $H(\sfu, t,x)$
is   continuous with respect to $\sfu$, is
  Lipschitz continuous with respect to $\sfu''$, and at 
all points of differentiability
of $H$ with respect to $\sfu''$ we have
$D_{\sfu''}H\in \bS_{\delta}$.

\end{assumption}

In the following theorem we need higher values of $p$
than in Theorem \ref{theorem 8,17.3} because in the proof
we need to use the embedding $W^{1,2}_{p}\subset
C^{0,1}$.

\begin{theorem}
                                     \label{theorem 8,18.2}
Let $p>d+2$ and suppose that Assumptions  \ref{assumption 8,18.3}
    and  \ref{assumption 7,16.5}  are satisfied and 
$\bar G\in L_{p}(\Pi)$.
  Then
there exist constants $ \theta,\hat\theta\in (0,1]$, depending   
only on $d$, $p$, $\delta$,  $K_{F}$, and $M_{2}(\Omega)$,
such that, if Assumptions \ref{assumption 8,18.2} and 
\ref{assumption 7,16.01}
are satisfied with these
 $\theta$ and $\hat\theta$, respectively,  then
there exists  
$u\in W^{1, 2}_{p }(\Pi ) $
 satisfying \eqref{8,17.7} in $\Pi $ \(a.e.\)  and such
that $u=g $ on $\partial'\Pi $. 
\end{theorem}

\begin{remark}
                                          \label{remark 8,20.1}
Observe that
generally there is no  uniqueness  
in Theorem \ref{theorem 8,18.2}.
For instance, in the one-dimensional case
the (quasilinear) equation 
$$
\partial_{t}u+ D^{2}u -(1-t)\sqrt{12|Du |}+2\sqrt{ (1-|x|^{3}) u} =0
$$
 in $\Pi=[0,1)\times (-1,1)$
with zero boundary data on $\partial'\Pi$ has two solutions:
one is identically equal to zero and the other one is
 $(1-t)^{2}(1-|x|^{3})$.

Uniqueness of solutions can be investigated
by using the results in \cite{Kr_85}.
\end{remark}

\begin{remark}
                                       \label{remark 2,21,1}
In case of linear equations Theorem \ref{theorem 8,18.2}
contains (apart from the restrictions on $p$)
the corresponding result of \cite{BC_93}
proved for equations with VMO main coefficients.

In Theorem 5.9 of Wang \cite{Wa92} 
one can find an a priori estimate 
for {\em any\/} viscosity solution
in case $H$
is independent of $\sfu'$ and $\Pi=C_{1}$.

By the way, it can be seen from our proofs that,
if $H$ is independent of $[\sfu']:
=(\sfu'_{1},...,\sfu'_{d})$, we can take $p>d+1$ in
Theorem \ref{theorem 8,18.2}.
\end{remark}

\begin{example}
For   $\tau>0$ take
$$
H(\sfu)=(1+\tau \cos\sqrt{|\ln |\sfu''||})\,Ê{\rm trace} \,\sfu'',
$$
and choose $\tau$ so small that $D_{\sfu''}H\in \bS_{\delta}$
for a $\delta\in(0,1]$. Then again $H$
is neither convex nor concave with respect to $\sfu''$
and our assumptions are satisfied perhaps with a further reduced
$\tau$ for 
$\bar F(\sfu'')={\rm trace} \,\sfu''$. An interesting
feature of this example is that, for generic $\sfu$,
 the limit of $(1/ \lambda )
H( \lambda \sfu)$ as $ \lambda \to\infty$ does not exist.

\end{example}

\begin{example}
                                             \label{example 1,13,1}
Let $A$ and $B$ be some countable sets and
assume that for $\alpha\in A$, $\beta\in B$,   $(t,x)\in\bR^{d+1}$,
  and $\sfu'\in\bR^{d+1}$
we are given an $\bS_{\delta}$-valued function  $a^{\alpha }(\sfu'_{0},t,x)$
(independent of $\beta$)
and a real-valued function $b^{\alpha\beta}(\sfu',t,x)$.
 Assume that 
  these functions are measurable
in $t,x$, $a^{\alpha }$ and $b^{\alpha\beta}$ are
 continuous with respect to $\sfu' $  
uniformly with respect to $\alpha,\beta,t,x$,
and  \smallskip
$$
\big|b^{\alpha\beta}(\sfu',t,x)\big|\leq K_{0}
|\sfu'|+\bar G(t,x),
\smallskip$$
where $\bar G \in L_{p}(\Pi)$, $p>d+2$.

Consider equation \eqref{8,17.7},
where
$$
H(\sfu,t,x):=
\infsup_{\beta\in B\,\,\alpha\in A}
\Big[\sum_{i,j=1}^{d}
a^{\alpha }_{ij}\big(\sfu'_{0},t,x\big)\sfu''_{ij}+
b^{\alpha\beta}(\sfu',t,x)\Big].
\medskip$$
Our measurability, boundedness, and countability assumptions guarantee that
$H$ is measurable in  $t,x$ and Lipschitz continuous in $\sfu''$. One can
also easily check that 
at all points of differentiability $D_{\sfu''}H\in\bS_{\delta}$.
Next assume that there is an $R_{0}\in(0,\infty)$ such that
for any   $z_{0}\in\Pi$, $r\in(0,R_{0}]$, and $\sfu'_{0}\in\bR$ one can find
$\bar{a}^{\alpha}\in\bS_{\delta}$ (independent of $t,x$) such that   
$$
 \dashint_{\hat C_{r}(z_{0})}\sup_{\alpha\in A} 
|a^{\alpha}(\sfu'_{0},z)-\bar{a}^{\alpha}|\,dz\leq\theta,
\quad\Big(\dashint_{\Gamma}h\,dz:=|\Gamma|^{-1}\int_{\Gamma}h\,dz\Big),
 $$
where $\theta$ is taken from Theorem \ref{theorem 8,18.2}.

Then we claim that the assertions   of 
Theorem \ref{theorem 8,18.2}  
hold true and estimate \eqref{8,18.1} holds with $\tau_{0}=0$.  

To prove the claim introduce
$$
F(\sfu'_{0},\sfu'',t,x)=\sup_{\alpha\in A}
 \sum_{i,j=1}^{d}
a^{\alpha }_{ij}(\sfu'_{0},t,x)\sfu''_{ij},\quad G=H-F.
 $$
Notice that Assumption \ref{assumption 8,18.2}  is satisfied 
with $\tau_{0}=0$ and  
$$
\bar{F}(\sfu''):=\sup_{\alpha\in A}
 \sum_{i,j=1}^{d}\bar
a^{\alpha\ }_{ij} \sfu''_{ij}
 $$
because these functions are convex, positive homogeneous
of degree one
with respect to $\sfu''$ and, for $|\sfu''|=1$,
$$
\dashint_{ \hat C_{r}(z_{0})}\big|
F(\sfu'_{0},\sfu'',z)-\bar{F}(\sfu'')\big|\,dz
\leq  \dashint_{\hat C_{r}(z_{0})}
\sup_{\alpha\in A} 
\Big|
\sum_{i,j=1}^{d}
\big[a^{\alpha }_{ij}\big(\sfu'_{0},z\big)-
\bar{a}^{\alpha}\big]\sfu''_{ij}\Big|\,dz
 $$
$$
\leq \dashint_{ \hat C_{r}(z_{0})}
\sup_{\alpha\in A} 
\big|a^{\alpha}\big(\sfu'_{0},z\big)-\bar{a}^{\alpha} \big|\,dz\leq\theta.
 $$

On can easily check that the remaining
item (iv) in Assumptions \ref{assumption 8,18.2} 
and Assumption \ref{assumption 7,16.01} (with $\hat\theta=0$)
are satisfied as well and this proves our claim. 
Thus Theorem \ref{theorem 8,18.2} is applicable.

As a result we have a solvability theorem for \eqref{8,17.7},
which covers (apart from the restriction on $p$), as $A$ and $B$ are 
singletons, the first
result about solvability of linear parabolic equations
with VMO coefficients obtained by Bramanti and Cerutti in \cite{BC_93}.
In this singleton case we also consider quasilinear equations.

\end{example}

 In the following theorem Assumption \ref{assumption 8,18.2}
is not used.

\begin{theorem}
                                       \label{theorem 2,2,2}  
Let $p>d+2$ and suppose that Assumptions \ref{assumption 7,16.1},
\ref{assumption 7,16.5}, and \ref{assumption 8,18.3} 
are satisfied, $\bar G\in L_{p}(\Pi)$, and \eqref{1,19,500} 
holds true. Let
$P(\sfu'')$ be a convex function on $\bS$ such that
at each point of its differentiability $D_{\sfu''}P\in\bS_{\delta'}$,
where $\delta'\in(0,\delta]$. Also assume that for any
$a\in \bS_{\delta}$ and $\sfu''\in\bS$ we have
$$
a^{ij}\sfu''_{ij}\leq  P(\sfu'')+K,
$$
where $K$ is a constant. Then  the equation
$$
\partial_{t}u+\max\big(H[u],P[u]\big)=0
$$
\(a.e.\) in $\Pi$ with boundary condition $u=g$
on $\partial'\Pi$ has a solution $u\in W^{1,2}_{p }(\Pi)$.
\end{theorem}
 
Proof. Introduce
$$
\hat H(\sfu,t,x)=\max\big(H(\sfu,t,x), P(\sfu'')\big),\quad
\hat F(\sfu'',t,x)=P(\sfu'')-P(0),\quad \hat G=\hat H-\hat F.
 $$
Obviously Assumptions \ref{assumption 8,18.2} and \ref{assumption 7,16.5},
are satisfied for $\hat H$, $\hat F$, and $\hat F$
in place of $H$, $F$, and $\bar F$, respectively, with a   $K_{F}$,
$\tau_{0}=\theta =0$, and $\delta'$ in place of $\delta$.
Finally, for any $\sfu,t,x$,   
$$
  \hat G(\sfu,t,x)=\max\big(H(\sfu,t,x)-P(\sfu'') 
+P(0) , P(0)\big)\geq P(0),
 $$
where for an $a\in \bS_{\delta}$   \smallskip
$$
H(\sfu,t,x)-P(\sfu'')=H(\sfu,t,x)-H(\sfu',0,t,x)-P(\sfu'')+H(\sfu',0,t,x)
\medskip$$
$$
=a^{ij}\sfu''_{ij}-P(\sfu'')+H(\sfu',0,t,x)\leq K+H(\sfu',0,t,x),
\medskip$$
which  together with \eqref{1,19,500} 
shows that Assumption \ref{assumption 7,16.01}
is also satisfied with $\hat\theta=0$ 
and $\bar G+K+ \big|P(0)\big|$ in place of $\bar G$.

Hence, Theorem \ref{theorem 8,18.2} is applicable   
and our theorem is proved.       \qed

\mysection{Interior estimates
of integral oscillations  of $D^{2}u$}
                                             \label{section 4,3,1}
Let
$F(\sfu'')$ be a convex function of $\sfu''\in\bS$
 (independent of $(t,x)$)
such that

(i) $F(0)=0$,

(ii) at all points of differentiability
of $F$ we have $ D_{\sfu''}F \in\bS_{\delta}$,
where $\delta\in(0,1]$ is a fixed number.

 The following theorem is a particular case
of the results in \cite{Li}.

\begin{theorem}
                                         \label{theorem 7,28.2}
There exists and $\bar\alpha=\bar\alpha(d,\delta)\in(0,1)$
such that
for any $\alpha\in(0,\bar\alpha]$ and
 $g\in C(\overline{\partial' C_{2}})$ there exists a unique
$v\in C(\bar{C}_{2})\cap 
C^{2+\alpha }_{\loc}(C_{2})$ satisfying
\begin{equation}
                                              \label{7,28.2}
\partial_{t}v+F(D^{2}v)=0\quad\text{in}\quad C_{2},\quad v=g\quad
\text{on}\quad\partial' C_{2}.
\end{equation}
Furthermore,
$$
|D^{2}v(z_{1})-D^{2}v(z_{2})|\leq N\rho^{\alpha }(z_{1},z_{2})
\sup_{\partial' C_{2}}|g|
$$
as long as $z_{1},z_{2}\in C_{1}$,
where   $N$ depends only on $\delta,\alpha$, and $d$.
\end{theorem}

Below in this section we fix
 $
\alpha\in(0,\bar\alpha]
 $. Recall that for a measurable set $\Gamma\subset\bR^{d+1}$
we denote by $|\Gamma|$ its Lebesgue measure,
and if $|\Gamma|\ne0$ and $u$ is integrable over $\Gamma$
we set
$$
u_{\Gamma}=\dashint_{\Gamma}u\,dxdt=\frac{1}{|\Gamma|}
\int_{\Gamma}u\,dxdt.
$$

\begin{lemma}
                                           \label{lemma 7,29.1}

Let $r\in(0,\infty)$, $\nu\geq2$ and let $\phi\in
C(\overline{\partial' C_{\nu r}})$. Then there exists
a unique $v\in C(\bar{C}_{\nu r})\cap C^{2+\alpha}_{\loc}
(C_{\nu r})$ such that
$$
\partial_{t}v+
F(D^{2}v)=0\quad\text{in}\quad C_{\nu r},\quad v=\phi\quad
\text{on}\quad\partial' C_{\nu r}.
$$
Furthermore,
$$
\dashint_{C_{r}}\dashint_{C_{r}}|D^{2}v(z_{1})-D^{2}v(z_{2})|
\,dz_{1}dz_{2}\leq N(d,\alpha,\delta)\nu^{-2-\alpha} r^{-2}\sup_{\partial'
C_{\nu r}}|\phi|.
$$
\end{lemma}

Proof. Scalings show that it suffices to concentrate on
$r=2/\nu$. In that case the existence of solution
follows from Theorem \ref{theorem 7,28.2}, which also implies that
for $z_{1},z_{2}\in C_{2/\nu}\subset C_{1}$
$$
|D^{2}v(z_{1})-D^{2}v(z_{2})|\leq N\nu^{-\alpha}\sup_{\partial'C_{2}}|\phi|.
$$
It only remains to observe that
$$
\dashint_{C_{2/\nu}}\dashint_{C_{2/\nu}}|D^{2}v(z_{1})-D^{2}v(z_{2})|
\,dz_{1}dz_{1}\leq\sup_{z_{1},z_{2}\in C_{2/\nu}}
|D^{2}v(z_{1})-D^{2}v(z_{2})|.
$$
  The lemma is proved.   \qed

Here is Theorem 1.9  of \cite{Kr_12.3}
combined with Theorem 2.3 of \cite{Kr_12.3}
(see also \cite{DKL_12}).

\begin{theorem}
                       \label{theorem 8.16.1}
Let  
$u\in   C(\bar C_{1})\cap  W^{1,2}_{d+1,loc}(C_{1})$. 
Then there are
constants $\bar\gamma=\bar\gamma(d,\delta,K)\in(0, 1]$
 and $N$,
depending only on
$\delta, d$, and $K$, such that for any $\gamma\in(0,\bar\gamma]$ and
  any operator $\cL=a^{ij}D_{ij}+b^{i}D_{i}$, with measurable
$\bS_{\delta}$-valued coefficients $a^{ij}$ 
and $b^{i}$, such that $|(b^{i})|\leq K $, given in $C_{1}$,
  we have
\begin{equation}
                                                          \label{8.11.01}
\int_{C_{1}} \big(|D^2u|^\gamma +|D u|^\gamma\big)\, dx \,  dt
 \leq  N \sup_{\partial'C_{1}}
 |u|^\gamma + N
\left(\int_{C_{1}}| \partial_{t}u+ \cL u|^{d+1} \,
dx\, dt\right)^{\gamma/{(d+1)}}.
\end{equation}  
\end{theorem}

Below we take
 $
\gamma\in(0,\bar\gamma]
 $.

\begin{lemma}
                                      \label{lemma 7,29.2}
 
Let $r\in(0,\infty)$ and $\nu\in[2,\infty)$. Then for any
$u\in W^{1,2}_{d+1 }(C_{\nu r}) $   we have  
$$
\big(\dashint_{C_{r} }\dashint_{C_{r}}
|D^{2}u(z_{1})-D^{2}u(z_{2})|^{\gamma}\,dz_{1}dz_{2}\big)^{1/\gamma}
$$
$$
\leq N\nu^{(d+2)/\gamma}
\big(\dashint_{C_{\nu r} }
|\partial_{t}u+F[u]|^{d+1}\,dz\big)^{1/(d+1)}
$$
\begin{equation}
                                            \label{7,29.3}
+
N\nu^{-\alpha}\big(\dashint_{C_{\nu r} }|D^{2}u|^{d+1}\,dz
\big)^{1/(d+1)},
\end{equation}
where $N$ depends only on $d,\alpha$, and $\delta$.
\end{lemma}

Proof. Define
$v$ to be a unique $C(\bar{C}_{\nu r})
\cap C^{2+\alpha}_{\loc}(C_{\nu r})$-solution of
equation $\partial_{t}v+F[v]=0$ in $C_{\nu r}$ with boundary condition
$v=u$ on $\partial' C_{\nu r}$. Such a function exists
by Lemma \ref{lemma 7,29.1}.
Furthermore, $v(x)-b^{i}x_{i}-c$ satisfies the same equation
for any constant $b^{i},c$. Hence by Lemma \ref{lemma 7,29.1}
and H\"older's inequality
$$
I_{r} :=\Big(
\dashint_{C_{r}}\dashint_{C_{r} }|D^{2}v(z_{1})-D^{2}v(z_{2})|^{\gamma}
\,dz_{1}dz_{2}\Big)^{1/\gamma}
$$
$$
\leq N \nu^{-2-\alpha} r^{-2}\sup_{z=(t,x)\in\partial'
C_{\nu r} }|u(z)-(D_{i}u)_{C_{\nu r}}x_{i}-u_{C_{\nu r}}|.
$$
By Poincar\'e's inequality (see, for instance,
Corollary 5.3 in \cite{DKL_12}) the last supremum is dominated
by a constant times
$$
\nu^{2} r^{2}\Big(\dashint_{C_{\nu r} }|D^{2}u|^{d+1}\,dz
\Big)^{1/(d+1)}.
$$
It follows that
\begin{equation}
                                                    \label{7,29.5}
I_{r} \leq N\nu^{-\alpha}\Big(\dashint_{C_{\nu r} }|D^{2}u|^{d+1}\,dz
\Big)^{1/(d+1)}.
\end{equation}

Next, the function $w:=u-v$ is of class $ 
 W^{1,2}_{d+1,\loc }(C_{\nu r})\cap C(\bar C_{\nu r})$
 and for an operator  
 $\cL=a^{ij}D_{ij}$ we have 
$$
\partial_{t}u+F[u]=\partial_{t}u+F[u]-(\partial_{t}v+F[v])
=\partial_{t}w+\cL w. 
$$
 Moreover, $w=0$ on $\partial' C_{\nu r}$.
Therefore,
by Theorem \ref{theorem 8.16.1},  there exists 
 $N=N(d,\delta)<\infty$ such that
$$
\dashint_{C_{r} }|D^{2}w|^{\gamma}\,dz
\leq \nu^{d+2}\dashint_{C_{\nu r} }|D^{2}w|^{\gamma}\,dz
$$
$$
\leq N\nu^{d+2}\Big(\dashint_{C_{\nu r} }
|\partial_{t}u+F[u]|^{d+1}\,dz\Big)^{\gamma/(d+1)}.
$$

Upon combining this result with \eqref{7,29.5} we come 
to \eqref{7,29.3} and the lemma is proved.   \qed

\mysection
{A priori estimates in $W^{1,2}_{p,\loc}$}
                                 \label{section 7,29.3}

Here we suppose that Assumptions  
\ref{assumption 7,16.01}  and \ref{assumption 8,18.2}
 are satisfied.
Thus, we assume that all assumptions on $H$ and $F$
stated before Theorem
\ref{theorem 8,17.3}
are satisfied. 
Take $
\alpha\in(0,\bar\alpha]$ and $\gamma\in(0,\bar\gamma]$.
First we note the following.

\begin{lemma}
                                          \label{lemma 8,15.1}
 For any $q\in[1,\infty)$ and $\mu>0$
there is a $\theta=\theta(d,\delta,K_{F},\mu,q)>0$ such that,
if Assumption \ref{assumption 8,18.2} is satisfied with this $\theta$,
then
for any $\sfu'_{0}\in\bR$, $r\in(0,R_{0}]$ and $z_{0}\in\Pi $ such that
$C_{r}(z_{0})\subset\Pi $
$$
 \dashint_{  C_{r}(z_{0})}\sup_{\substack{\sfu''\in\bS,\\|\sfu''|>\tau_{0}}}
\frac{\big|F(\sfu'_{0},\sfu'',z)-\bar{F}(\sfu'')\big|^{q}}
{|\sfu''|^{q}}\,dz\leq \mu^{q},
$$
where $\bar F=\bar F_{z_{0},r,\sfu'_{0}}$.
\end{lemma}

The proof of this lemma is practically identical to that
of Lemma 5.1 of \cite{Kr_13}
given there for the elliptic case. 

\begin{lemma}
                                          \label{lemma 11.10.1}
Let $u\in W^{1, 2}_{d+1,\loc}(\Pi)$. Then there exist an
$\bS_{\delta}$-valued function $a(t,x)$,   $\bR^{d}$-valued
functions $b(t,x)$, and real-valued function  
  $f(t,x)$, 
such that they are measurable,
$$
|b|  \leq K_{0},  \quad| f|\leq\bar G+K_{0}|u|,
$$ 
and in $\Pi$ \(a.e.\)
\begin{equation}
                                              \label{11.10.3}
a^{ij}D_{ij}u+  b^{i}D_{i}u +  f=H[u].
\end{equation}
\end{lemma}

This is a simple consequence of the fact that
 there is an $\bS_{\delta}$-valued function $a $
such that
$$
H[u](t,x)-H(u,Du,0,t,x)=a^{ij}D_{ij}u,
$$
and
$$
|H(u,Du,0,t,x)|\leq K_{0}(|u|+ |Du|)+\bar G.
$$

\begin{lemma}
                                       \label{lemma 7,29.3}
Let $r\in(0,\infty)$ and
 $\nu\geq 2$ satisfy  $\nu r\leq 
 R_{0} $.
Take   
$$
\mu\in(0,\infty),\quad \beta\in(1,\infty),
$$ 
and suppose that Assumption \ref{assumption 8,18.2}
is satisfied with $\theta=\theta(d,\delta,K_{F},\mu,\beta d+\beta)$
\(see Lemma \ref{lemma 8,15.1}\).
Take a function $u\in W^{1,2}_{d+1}(\Pi )$  and 
for $z_{0}\in\Pi $
 such that $C_{\nu r }(z_{0}) \subset\Pi $ \(if such $z_{0}$'s
exist\) denote
$$
I_{r}(z_{0})=\big(\dashint_{C_{r}(z_{0}) }
\dashint_{C _{r}(z_{0}) }
|D^{2}u(z_{1})-D^{2}u(z_{2})|^{\gamma}\,dz_{1}dz_{2}\big)^{1/\gamma}.
$$

Then   
$$
I_{r} (z_{0})\leq  N\nu^{(d+2)/\gamma}\bigg(
\dashint_{C_{\nu r }(z_{0}) } \big|\partial_{t}u+F[u]\big|^{d+1}
 \,dz\bigg)^{1/(d+1)}+N\tau_{0}\nu^{(d+2)/\gamma}
$$

\begin{equation}
                                              \label{7,29.6}
+N\Big[\big(\mu+\omega_{F,u,\Pi}(\nu r )\big)\nu^{(d+2)/\gamma} +\nu^{-\alpha} 
 \Big]
\bigg(
\dashint_{C_{\nu r }(z_{0}) }
|D^{2}u |^{\beta' (d+1)}\,dz\bigg)^{1/(\beta'(d+1))},
\end{equation}
where $\beta'=\beta/(\beta-1)$  and 
$N$ depends only on $ d,K_{F}$, $\alpha$,
and $\delta$.

\end{lemma}

Proof.   Set $\rho:=\nu r$. Since    $\rho\leq R_{0}$,
 $\bar{F}=\bar{F}_{z_{0},\rho,u(z_{0})}$ is well defined
 and  by Lemma \ref{lemma 7,29.2}
$$
I_{r}(z_{0})
\leq N\nu^{(d+2)/\gamma}\bigg(\dashint_{C_{\rho}(z_{0}) }
 \big|\partial_{t}u+\bar{F}[ u]\big|^{d+1} \,dz\bigg)^{1/(d+1)}
$$
\begin{equation}
                                                     \label{11.27.7}
+N \nu^{-\alpha} 
\bigg(\dashint_{C_{\rho}(z_{0}) }
|D^{2}u |^{d+1}\,dz\bigg)^{1/(d+1)}.
\end{equation}

 By setting $\hat F[u](z)=F(u(z_{0}),D^{2}u(z))$ we find 
$$
 \dashint_{C_{\rho}(z_{0}) }
\big|\partial_{t}u+\bar{F}[u]\big|^{d+1}\,dz \leq
N \dashint_{C_{\rho}(z_{0}) }
\big|\partial_{t}u+F[u]\big|^{d+1}\,dy + NJ_{1}+NJ_{2},
$$
where 
$$
J_{1}=\dashint_{C_{\rho}(z_{0}) }
\big|\hat F[u]-\bar{F}[u]\big|^{d+1}\,dz
$$
is dominated by
$$
\dashint_{C_{\rho}(z_{0})  }I_{|D^{2}u|>\tau_{0}}\frac{
\big|\hat F[u]-\bar{F}[u]\big|^{d+1}}{|D^{2}u|^{d+1}}
|D^{2}u|^{d+1}\,dz+N\tau_{0}^{d+1},
$$
which in turn owing to Lemma \ref{lemma 8,15.1}
and H\"older's inequality is less than
$$
N\mu^{d+1}\bigg(\dashint_{C_{\rho}(z_{0}) }
|D^{2}u|^{\beta'(d+1)}\,dz\bigg)^{1/\beta'}+N\tau_{0}^{d+1},
$$
and 
$$
J_{2}=\dashint_{C_{\rho}(z_{0}) }
\big|\hat F[u]- F[u]\big|^{d+1}\,dz
\leq\omega^{d+1}_{F}(\osc_{C_{\rho}(z_{0}) }
u)\dashint_{C_{\rho}(z_{0}) }
|D^{2}u|^{d+1}\,dz.
$$

It follows that
$$
\bigg(\dashint_{C_{\rho}(z_{0}) }
\big|\partial_{t}u+\bar{F}[u]\big|^{d+1}\,dy\bigg)^{1/(d+1)}\leq
N \bigg(\dashint_{C_{\rho}(z_{0}) }
\big|\partial_{t}u+F[u]\big|^{d+1}\,dy \bigg)^{1/(d+1)} 
$$
$$
+N\mu \bigg(\dashint_{C_{\rho}(z_{0}) }
|D^{2}u|^{\beta'(d+1)}\,dy\bigg)^{1/(\beta'd+\beta')}+N\tau_{0} 
$$
$$
+N\omega_{F,u,\Pi}(\rho )\bigg(\dashint_{C_{\rho}(z_{0}) }
|D^{2}u |^{d+1}\,dz\bigg)^{1/(d+1)}.
$$
 
This and \eqref{11.27.7} yield \eqref{7,29.6} 
since 
$$
\bigg(\dashint_{C_{\rho}(z_{0}) }
|D^{2}u |^{d+1}\,dz\bigg)^{1/(d+1)}\leq
\bigg(\dashint_{C_{\rho}(z_{0}) }
|D^{2}u |^{\beta'(d+1)}\,dz\bigg)^{1/(\beta'(d+1))}
$$
 by H\"older's inequality. The lemma is proved.                    \qed

\begin{lemma}
                                        \label{lemma 7,30.1}
Take    $p>d+1$, $R\in(0,1]$,   and  $u\in W^{1,2}_{p}(C_{2R })$.
 Then
there exist    constants $\hat\theta,\theta\in (0,1]$, depending
only on $d$, $p$, $\delta$, and  $K_{F}$,
such that, if Assumptions \ref{assumption 7,16.01} and
 \ref{assumption 8,18.2} are satisfied with these $\hat\theta$ 
and $\theta$, respectively,
  then 
there is a constant  $N$,
depending only on $R_{0}$, $d$, $p$, $K_{0}$, $K_{F}$,  $\delta$, 
and    $\omega_{F,u,C_{2R}}$, such that  
$$
\|D^{2}u\|_{L_{p}(C_{R })}\leq
N \big\|\partial_{t}u+H[u]\big\|_{L_{p}(C_{2R })}+N\| \bar G\|_{L_{p}(C_{2R })}
+N\tau_{0}   
$$
\begin{equation}
                                                      \label{7,31.1}
+N  R^{(d+2)(1/p-1/\gamma) }\big\|
\,|D^{2}u|^{\gamma}
\big\|^{1/\gamma}_{L_{1}(C_{2R })}+N\| u\|_{L_{p}(C_{2R })},
\end{equation}

$$
\|D^{2}u\|_{L_{p}(C_{R})}\leq N\tau_{0}  +NR^{(d+2)/p-2}\sup_{C_{2R}}|u|
$$ 
\begin{equation}
                                              \label{7,31.2}
+N \big( \big\|\partial_{t}u+H[u]\big\|_{L_{p}(C_{2R})}+\|
 \bar G\|_{L_{p}(C_{2R})}\big).
\end{equation}

\end{lemma}

 Proof. 
  For   $\rho>0$, and $z \in
Q:=\bR_{+}\times\bR^{d}$ introduce
$$
I_{r}(h,z)=
\bigg(
\dashint_{C_{r}(z )}\dashint_{C_{r}(z  )}
|h(z_{1})-h(z_{2})|^{\gamma}\,dz_{1}dz_{2}\bigg)^{1/\gamma},
$$
$$
  h^{\shharp}_{Q,\gamma,\rho}(z)=\sup\{I_{r}( h,z_{0}):
z_{0}\in
Q,r\in(0,\rho],C_{r}(z_{0})\ni z\},
$$
\begin{equation}
                                            \label{8,12.6}
\bM  h(z ) =\sup_{\substack{r>0,\\ 
   C_{r}(z_{0})\ni z}}\dashint_{C_{r}(z_{0}) }
|h(\zeta)|\,d\zeta,
\end{equation}
\index{$h^{\shharp}_{Q, \gamma,\rho}$}%
\index{$\bM  h$}%
whenever these definitions make sense.
Note that 
 $h^{\shharp}_{Q,\gamma,\rho}$ is well defined
in $C_{R}$
for measurable $h$ even defined only in $C_{R+2\rho}$.

Then take $\varepsilon\in(0,1]$ to be specified later and
 take $R_{1}<R_{2}\leq 2R$ such that 
\begin{equation}
                                                 \label{7,31.3}
 R_{2}-R_{1}\leq \varepsilon  R_{0},\quad R_{2}\leq  2R_{1}.
\end{equation}
Next,
take $\nu\geq2$ and set 
$$
 r_{0}=(R _{2}-R_{1})/(\nu+1) .
$$
Observe that $\nu r_{0}\leq \varepsilon  R_{0}$ and $R_{2} -\nu r_{0}=R_{1}+r_{0}$.
It follows that,   if $r\leq r_{0}$, $z\in C_{R_{1}}$,
and $z\in C_{r}(z_{0})$, then $C_{\nu r}(z_{0})\subset
C_{R_{2}}$,
 which by Lemma \ref{lemma 7,29.3} applied with $\Pi =C_{R_{2}}$
implies that  
$$
I_{r}(z_{0})\leq N\nu^{(d+2)/\gamma}\bM^{1/(d+1)}
\big(\big|\partial_{t}u+F[u]\big|^{d+1}I_{
C_{R_{2} }}\big)(z)+N\tau_{0}\nu^{(d+2)/\gamma}
$$

$$
+N\Big[\big(\mu+\omega_{F,u, C_{2R}}( \nu r_{0})\big)\nu^{(d+2)/\gamma} +\nu^{-\alpha}\Big]
\bM^{1/(\beta'(d+1))}\big(|D^{2}u|^{\beta'(d+1)}I_{
C_{R_{2} }}\big)(z)
$$
with $N$ depending only on   $ d,K_{F}$, 
and $\delta$.
It follows that in $C_{R_{1}}$
$$
(D^{2}u)^{\shharp}_{Q,\gamma,r_{0}} \leq N\nu^{(d+2)/\gamma}
\bM^{1/(d+1)}\big(\big|\partial_{t}u+F[u]\big|^{d+1}I_{
C_{R_{2} }}\big) +N\tau_{0}\nu^{(d+2)/\gamma}
$$

$$
+N\Big[\big(\mu+\omega_{F,u,C_{2R}}(\varepsilon R_{0})\big)\nu^{(d+2)/\gamma}
 +\nu^{-\alpha} 
\Big]
\bM^{1/(\beta'(d+1))}\big(|D^{2}u|^{\beta'(d+1)}I_{
C_{R_{2} }}\big) .
$$
By Theorem \ref{theorem 8,4.1},   
 with 
$$
\kappa=r_{0}/R_{1} \leq
1/3 ,\quad \chi_{1}=(d+2)/\gamma, \quad\chi_{2}=
(d+2)(1/\gamma-1/p )
$$
and the Hardy-Littlewood maximal function theorem, by taking $\beta$ so that
  $p>\beta'(d+1)$, we obtain
$$
\|D^{2}u\|_{L_{p}(C_{R_{1}})}\leq
N \nu^{(d+2)/\gamma}\big\|F[u]\big\|_{L_{p}(C_{R_{2} })}
+N\tau_{0}\nu^{(d+2)/\gamma}|R_{1}|^{(d+2)/p}
$$
$$
+\Big[N \big(\mu+\omega_{F,u,C_{2R}}(\varepsilon R_{0})\big)
\nu^{(d+2)/\gamma}+N_{0}\nu^{-\alpha}\Big] 
\|D^{2}u\|_{L_{p}(C_{R_{2} })}
$$
\begin{equation}
                                          \label{4,21,3}
+N\nu^{\chi_{1}}(R_{2}-R_{1})^{-\chi_{1}}R_{1}^{-\chi_{2}+\chi_{1}}
\big\|\,|D^{2}u|^{\gamma}
\big\|^{1/\gamma}_{L_{1}(C_{2R })},
\end{equation}
where and below the constants $N$, $N_{i}$ depend only on $ d$, $p$, $K_{F}$, 
and $\delta$.  

Now we take and fix $\nu\geq2$ so that
$$
N_{0}\nu^{-\alpha}\leq 1/4.
$$
Then \eqref{4,21,3} becomes
$$
\|D^{2}u\|_{L_{p}(C_{R_{1}})}\leq
N_{1}  \big\|F[u]\big\|_{L_{p}(C_{R_{2} })}
+N\tau_{0} |R_{1}|^{(d+2)/p}
$$
$$
+\Big[N_{2} \big(\mu+\omega_{F,u,C_{2R}}(\varepsilon R_{0})\big)
 +1/4\Big] 
\|D^{2}u\|_{L_{p}(C_{R_{2} })}
$$
\begin{equation}
                                          \label{4,21,4}
+N\nu^{\chi_{1}}(R_{2}-R_{1})^{-\chi_{1}}R_{1}^{-\chi_{2}+\chi_{1}}
\big\|\,|D^{2}u|^{\gamma}
\big\|^{1/\gamma}_{L_{1}(C_{2R })},
\end{equation}

Next, we use the fact that   \smallskip
$$
\big|F[u]\big|\leq\big|H[u]\big|+K_{0}|u|+K_{0}|Du| +\bar G+\hat\theta|D^{2}u|
\smallskip$$
and that by interpolation inequalities  \smallskip
$$
K_{0}N_{1} \|D u\|_{L_{p}(B_{R_{2} })}\leq 
(1/8)
\|D^{2}u\|_{L_{p}(B_{R_{2} })}
+
N  
\|u\|_{L_{p}(B_{R_{2} })}.
\smallskip$$
Then we take $\hat \theta$ and $\mu$ so small that
$$
N_{1}\hat\theta\leq 1/8,\quad N_{2}\mu\leq 1/8,
$$
and, finally, take 
the largest $\varepsilon\leq1$ such 
that   \smallskip
$$
N_{2}
\omega_{F,u,C_{2R}}(\varepsilon R_{0})\leq 1/8.
\smallskip$$
This $\varepsilon$ will appear later
in our arguments and this is the way how the constant $N$ in the statement
of the lemma depends on $\omega_{F,u,C_{2R}}$.

Then we require that Assumptions \ref{assumption 7,16.01}
 and \ref{assumption 8,18.2}
be satisfied with the above chosen $\hat\theta$ and
 $\theta=\theta(d,\delta,K_{F},\mu,\beta d+\beta)$
 (see Lemma \ref{lemma 8,15.1}), respectively. By combining the above
  we get
$$
\|D^{2}u\|_{L_{p}(C_{R_{1}})}\leq
N \big\|\partial_{t}u+H[u]\big\|_{L_{p}(C_{R _{2}})}
+N\tau_{0}   R  ^{d/p} 
+(5/8)
\|D^{2}u\|_{L_{p}(C_{R_{2} })}
$$
\begin{equation}
                                                    \label{8,20.3}
+N(R_{2}-R_{1})^{-\chi_{1}}R_{1}^{-\chi_{2}+\chi_{1}}
\big\|\,|D^{2}u|^{\gamma}
\big\|^{1/\gamma}_{L_{1}(C_{2R })}+N\|u\|_{L_{p}(C_{2R  })}
+N\| \bar G\|_{L_{p}(C_{2R })}.
\end{equation}

 Now we are going to iterate this
estimate by defining $R_{1}=R$ and for $k\geq 1$  
$$
R_{k+1}=R_{k }+c R (n_{0}+k)^{-2},
$$
where the constant $c=O(n_{0})$ is chosen so that $R_{k}
\uparrow 2R$ as $k\to\infty$, that is
$$
c\sum_{k=1}^{\infty}(n_{0}+k)^{-2}=1.
$$
and 
$n_{0}$ is chosen so that for $k\geq 1$
$$
R_{k+1}-R_{k}=c R (n_{0}+k)^{-2}\leq  R cn^{-2}_{0}
\leq R\leq  R_{k},
$$
which is satisfied if $n_{0}$ is just an absolute
constant, and  (this time we need $n^{-1}_{0}=o(
\varepsilon R_{0} )$ as $\varepsilon  R_{0}\to0$)
$$
R_{k+1}-R_{k}=cR(n_{0}+k)^{-2}\leq  cn^{-2}_{0}\leq \varepsilon R_{0}.
$$
Also observe that $R\leq R_{k}\leq 2R$ and
$$
(R_{k+1}-R_{k})^{-\chi_{1}} R_{k}^{-\chi_{2}+\chi_{1}}
\leq N(n_{0}+k)^{2\chi_{1}} R ^{-\chi_{2}} .
$$

Then for $k\geq1$ we get
$$
\|D^{2}u\|_{L_{p}(C_{R_{k}})}\leq
N \big\|\partial_{t}u+H[u]\big\|_{L_{p}(C_{R _{k+1}})}
+N\tau_{0}   R  ^{d/p} 
+(5/8)
\|D^{2}u\|_{L_{p}(C_{R_{k+1} })}
$$
$$
+N (n_{0}+k)^{2\chi_{1}}R^{-\chi_{2}}\big\|
\,|D^{2}u|^{\gamma}
\big\|^{1/\gamma}_{L_{1}(C_{2R })}+N\| u\|_{L_{p}(C_{2R })}
+N\| \bar G\|_{L_{p}(C_{2R })}.
$$
We multiply both parts of this inequality by $(5/8)^{ k}$
and sum up the results over $k=1,2,...$. Then we cancel
like terms
$$
\sum_{k=2}^{\infty}(5/8)^{ k}\|D^{2}u\|_{L_{p}(C_{R_{k}})},
$$
which are finite since $u\in W^{1,2}_{p}(B_{2R})$,
and finally take into account that
$$
 \sum_{k=2}^{\infty}(5/8)^{ k}(n_{0}+k)^{2\chi_{1}}
\leq N n_{0}^{2\chi_{1}}\sum_{k=2}^{\infty}(5/8)^{ k}
+N\sum_{k=2}^{\infty}(5/8)^{ k}k ^{2\chi_{1}}\leq N.
$$
Then we come to \eqref{7,31.1}.

Next, by using equation  \eqref{11.10.3}    
and performing scaling 
in Theorem \ref{theorem 8.16.1} (here we need $R\leq1$),
 using H\"older's
inequality (to go from $d+1$ to $p$),
and denoting
$$
I=\|\bar G \|_{L_{p}(C_{2R})}+
\big\|\partial_{t}u+H[u] \big\|_{L_{p}(C_{2R})}
$$
 we infer that in \eqref{7,31.1}
$$
\big\|\,|D^{2}u|^{\gamma}
\big\|^{1/\gamma}_{L_{1}(C_{2R})}
\leq NR^{\chi_{2}}
\Big(\|\bar G+K_{0}|u|\,\|_{L_{p}(C_{2R})}+
\big\|\partial_{t}u+H[u] \big\|_{L_{p}(C_{2R})}\Big)
$$
$$
+NR^{\chi_{3}}\sup_{C_{2R}}|u|\leq N R^{\chi_{2}}
I+NR^{\chi_{3}}\sup_{C_{2R}}|u|,
$$
where $\chi_{3}=(d+2)/\gamma-2$.
After that it suffices to   roughly estimate 
$\| u\|_{L_{p}(C_{2R})}$ in \eqref{7,31.1} by the last term above. 
The lemma is proved.                   \qed

\mysection 
{Boundary a priori estimates in the simplest case}
                                    \label{section 8.19.2}

Introduce
$$
\bR^{d+1}_{+}=\{(t,x):t\in\bR,x=(x^{1},...,x^{d})\in\bR^{d},x^{1}>0\},
$$
$$
B_{r}^{+}(x_{0})=B_{r}(x_{0})\cap\{x^{1}>0\},\quad
C^{+}_{\tau,r}(t_{0},x_{0})= [t_{0},t_{0}+\tau)\times 
B^{+}_{r}(x_{0}), 
$$
$$
\partial_{x^{1}}C^{+}_{\tau,r}(t_{0},x_{0})=
\bar C^{+}_{\tau,r}(t_{0},x_{0})\cap\{ x^{1}=0\},
$$
where $\tau,r\geq0$,   $(t_{0},x_{0})\in\bar{\bR}^{d+1}_{+}$. If $t_{0}=0,x_{0}=0$,
 we drop $(t_{0},x_{0})$ in the arguments above. Also, 
if $\tau=r^{2}$ we write
$r$ in place of $\tau,r$ in the subscripts, for 
instance,
$$
C^{+}_{ r} (t_{0},x_{0}) :=C^{+}_{r^{2},r}(t_{0},x_{0}).
$$

Take $ \gamma  $  from Section \ref{section 7,29.3}
and $\alpha\in(0,1)$ to be determined later.
Let $F$ be the function from Section \ref{section 4,3,1}.

\begin{lemma}
                                             \label{lemma 8,17.3}
If
 $  r>0$, $z_{0} \in\bar\bR^{d+1}_{+}$,
$\nu\geq  12 $,
$$
u\in  \bigcap_{\rho<\nu r}
W^{1,2}_{d+1}(C^{+}_{ \rho}(z_{0}))\cap 
C(\bar C_{ \nu r}^{+}(z_{0})),
$$
and $u $ vanishes on 
$\partial_{x^{1} }C_{ \nu r}^{+}(z_{0})$
if this set is nonempty,  then
we have
$$
 \Big(\dashint_{C^{+}_{ r}(z_{0})}\dashint_{C^{+}_{ r}(z_{0})}
|D^{2}u(z_{1})-D^{2}u(z_{2})|^{\gamma}
\,dz_{1}dz_{2}\Big)^{1/\gamma}
$$
$$
\leq N\nu^{(d+2)/\gamma}
\Big(\dashint_{C^{+}_{ \nu r}(z_{0})}
\big|\partial_{t}u+F[u]\big|^{d+1}\,dz\Big)^{1/(d+1)}
$$
\begin{equation}
                                              \label{8,17.3}
+
N\nu^{-\alpha}\Big
(\dashint_{C^{+}_{ \nu r}(z_{0})}|D^{2}u|^{d+1}\,dz
\Big)^{1/(d+1)},
\end{equation}
where $N$ depends only on $d$ and $\delta$.

\end{lemma}

Proof.
 Scalings show that it suffices to prove
the lemma only for $\nu r=3$. 
Furthermore, without loss of generality
we may assume that $z_{0}=(0,x_{0})$ and
$x_{0}=(|z_{0}|,0,...,0)\in\bR^{d}$.
Then
we consider two cases.  

{\em Case 1: $|z_{0}|>1/2$.} In this case, we have 
$$
B^+_{ r }(x_{0})
= B^{+}_{3/\nu}(x_{0})= B _{r}(x_{0})\subset
 B_{ \nu'r}(x_{0}) \subset \bR^d_+,\quad
C_{\nu'r}(z_{0})\subset \bR^{d+1}_+,
$$
where $\nu'=\nu/6$ ($\geq 2$).
 Therefore, inequality
\eqref{8,17.3} is an immediate consequence of Lemma 
\ref{lemma 7,29.2}.

{\em Case 2: $|z_{0}|\in [0,1/2]$.} Since $r=3/\nu\leq1/2$, we have
$$
B^{+}_r(x_{0})\subset B^{+}_1 \subset B^{+}_2\subset
B^{+}_{3}(x_{0})= B^{+}_{\nu r}(x_{0}).
$$

Let $v$ be the classical solution of 
$\partial_{t}v+F[v]=0$ in $C^+_2 $
with boundary condition $v=  u$ on $\partial' C _{2}^{+}
 $.
Such a solution exists due to the results in \cite{Li},
which also provide  an estimate on $D^{2}v$, so that
(for $\alpha\in(0,\alpha_{0}(d,\delta)]$)
$$
I:=\dashint_{C^{+}_r(z_{0})}\dashint_{C^{+}_r(z_{0})}|
  D^2v (z_{1})-
  D^2 v (z_{2})|\,dz_{1}dz_{2} 
\le Nr^\alpha[  D^2v]_{ C^{\alpha}(C^{+}_1)}
$$
$$
\le Nr^\alpha\sup_{C^+_2} |v|
 =Nr^\alpha\sup_{\partial' C^+_2} | u|
$$
where the last equality is a consequence of the maximum
principle and the fact that $F(0)=0$.
By employing
 Poincar\`e's inequality ($u=0$ on $\partial_{x^{1}}
C _{2}^{+}$), we see that 
$$
I\leq Nr^{\alpha}
\Big(\dashint_{B^{+}_2}\big
(|\partial_{t}u|^{d+1}+ |D^2 u|^{d+1}\big) \,dz\Big)^{1/(d+1)}.
$$
Here $r^{\alpha}=N\nu^{-\alpha}$ and
$$
|\partial_{t}u|\leq\big|\partial_{t}u+F[u]\big|+\big|F[u]\big| 
\leq \big|\partial_{t}u+F[u]\big|+N|D^{2}u|.
$$
Therefore,
$$
I\leq N\nu^{-\alpha}\Big(\dashint_{B^{+}_2} 
 |D^2 u|^{d+1} \,dz\Big)^{1/(d+1)}
+N\Big(\dashint_{C^{+}_{ 2}}
\big|\partial_{t}u+F[u]\big|^{d+1}\,dz\Big)^{1/(d+1)}
$$

Next,
recall that $\gamma\in (0,1]$. By H\"older's inequality, we get
$$
\dashint_{C^{+}_r(z_{0})}\dashint_{C^{+}_r(z_{0})}|
   D^2v (z_{1})-   D^2 v (z_{2})|^\gamma\,dz_{1}dz_{2}
$$
$$
\le N \nu^{-\gamma\alpha}
\Big(\dashint_{C^{+}_{\nu r}(z_{0})} | D^2 u|^{d+1}
 \,dz\Big)^{\gamma/(d+1)}
$$
\begin{equation}
                                            \label{8,17.1}
\phantom{nnnnn}\quad +N\Big(\dashint_{C^{+}_{\nu r}(z_{0})}
\big|\partial_{t}u+F[u]\big|^{d+1}\,dz\Big)^{\gamma/(d+1)}.
\end{equation}

Next, use again that 
$$
f:=\partial u+F(D^{2} u)=\partial(u-v)+F(D^{2} u)-F(D^{2}v)
=\partial w+a^{ij}D^{2}_{ij}w  
$$
in $C^{+}_{2}$ and $w=0$ on $\partial' C^{+}_{2}$,
where $(a^{ij})$ is an $\bS_{\delta}$-valued function and 
$w=  u-v$.
We extend $f$ and $w$ to all of $C_{2}$ as  odd functions of $x^{1}$
and adjust $a^{ij}$ appropriately so as to have
equation $f=\partial w+a^{ij}D^{2}_{ij}w $ in $C_{2}$, to which we
apply Theorem \ref{theorem 8.16.1} and get (recall that $\nu r=3$)
$$
\dashint_{C^{+}_r(z_{0})}|  D^2w|^{\gamma}\,dz
\leq N r^{-d-2}\int_{C^{+}_2}|  D^2w|^\gamma\,dz
$$
$$
\leq N \nu^{d+2}\Big(\dashint_{C^{+}_{\nu r}(z_{0})}
 |f|^{d+1}\,dz\Big)^{\gamma/(d+1)}
$$
and
$$
\dashint_{C^{+}_r(z_{0})}
\dashint_{C^{+}_r(z_{0})}|  D^2w (z_{1})-  D^2w (z_{2}
)|^{\gamma}\,dz_{1}dz_{2}  
$$
$$
\leq N \nu^{d+2}\Big(\dashint_{C^{+}_{\nu r}(z_{0})}
|f|^{d+1}\,dz\Big)^{\gamma/(d+1)}.
$$
Combining this with \eqref{8,17.1}
and observing that $D^{2}u=D^{2}v+D^{2}w$ yield \eqref{8,17.3}
in Case 2 as well. 
The lemma
is proved.   \qed

Coming back to our domain $\Omega$ recall that
we say that $\Omega$
is a  $C^{1,1}$-domain if  there exists $\rho_{0}=\rho_{0}(\Omega)
\in(0,1]$
 for which at any point
$x_{0}\in\partial\Omega$ there is an orthonormal system
of coordinates $\Psi(x_{0})$ with the origin at $x_{0}$
 such that in the new coordinates $\tilde{x}
=(\tilde{x}^{1},\tilde{x}')$  there exists a function
$$
\psi \in C^{1,1}(\{\tilde{x}'\in\bR^{d-1}:
|\tilde{x}'|\leq8\rho_{0}\})
$$
with the $C^{1,1}(B_{8\rho_{0}})$-norm majorated by a constant $M_{2}(\Omega)$ independent
of $x_{0}$ and such that
$$
\psi(0)=0,\quad \psi_{\tilde{x}_{i}}(0)=0,
\quad i=2,...,d,\quad |D_{x'}\psi(
\tilde x')|\leq1\quad\text{for}
\quad|\tilde x'|\leq 8\rho_{0},
$$
$$
 \{\tilde{x}:|\tilde{x}'|\leq 8\rho_{0},\,
\psi(\tilde{x}')+8\rho_{0}\leq\tilde{x}^{1}\leq\psi(\tilde{x}')+8\rho_{0}\}
\cap \Omega 
$$
$$
= \{\tilde{x}:|\tilde{x}'|\leq8\rho_{0},\,
\psi(\tilde{x}')<\tilde{x}^{1}\leq\psi(\tilde{x}')+8\rho_{0}\}.
$$

Below in this section we assume that
$$
0\in\partial\Omega
$$
and that the original system of coordinates in
$\bR^{d}$ coincides with the one described above
for $x_{0}=0$.
 
\begin{lemma}
                                   \label{lemma 2.14.1}
Introduce
$$
\Gamma :=\{x:|x'|\leq 8 \rho_{0}(\Omega) ,
\psi(x')<x^{1}\leq\psi(x')+8\rho_{0}(\Omega) \}\quad(
\subset \Omega ),
$$
$$
\hat{\Gamma} :=
\{y:|y'|\leq 8\rho_{0}(\Omega) ,0< y^{1}\leq
8\rho_{0}(\Omega) \}.
$$
Also introduce a mapping $x\to y(x)$ of $\Gamma  $ onto
$\hat{\Gamma} $ by
\begin{equation}
                                           \label{11.11.01}
x^{1}\to y^{1}=y^{1}(x)=x^{1}-\psi(x'),\quad
 x'\to y'=y'(x)=x'.
\end{equation}
Then this mapping  has an  inverse $y\to x(y)$. Furthermore,
the Jacobians of both mappings are equal to one. 
\end{lemma}

This lemma is obvious.
 
It is convenient to extend $\psi(x')$
for $|x'|\geq 8\rho_{0}(\Omega)$, so that the extension
is smooth and has the magnitude of the gradient
bounded by one and define $y(x)$ by the same formula 
\eqref{11.11.01} for all $x\in\bR^{d}$. Of course,
by $x(y)$ we mean the inverse of $y(x)$. 
Obviously, the assertions of Lemma \ref{lemma 2.14.1}
hold true
for such extensions.
\begin{remark}
                                     \label{remark 8,12.1}
For $r\in(0,\infty)$ and $z\in\bR^{d}$ define
$$
\chB^{+}_{r}=x(B^{+}_{r}),
\quad \chB^{+}_{r}(z)=x(B_{r}^{+}(y(z))) .
$$
Then, as is easy to see

(i) $\chB^{+}_{r}\subset\Gamma\subset\Omega$ if
$r\leq  8 \rho_{0}(\Omega)$;

(ii) $\chB^{+}_{r}(z)\subset \chB^{+}_{4\rho_{0}(\Omega)}$
if $\rho>0$,
$\rho+ r\leq  8 \rho_{0}(\Omega)$,
and $z\in\chB^{+}_{\rho}$ .

\end{remark}

\begin{lemma}
                                 \label{lemma 6.2.2}
Take $z\in\chB^{+}_{2\rho_{0}(\Omega)}$. Then

\(i\,\) for $r\leq 2\rho_{0}(\Omega)$   we have
\begin{equation}
                                   \label{6.2.4}
\chB^{+}_{r }(z)\subset B_{2r}(z)\cap\Omega,
\quad B_{r/2}(z)\cap\Omega\subset\chB^{+}_{r}(z);
\end{equation}

\(ii\,\) if $\nu\geq 1$  and $\nu r\leq 2\rho_{0}(\Omega)$, we have
\begin{equation}
                                         \label{6.2.5}
|\chB^{+}_{\nu r}(z)|\leq N(d)\nu^{d}
|\chB^{+}_{r}(z)|.
\end{equation}
 
\end{lemma}

Proof. (i). First notice that $\chB^{+}_{r }(z)\subset \Gamma$.
Then, since $|D_{x'}\psi|\leq1$,
for any $x_{1},x_{2}\in\Gamma$ we have $|y(x_{1})-
y(x_{2}|\leq 2|x_{1}-x_{2}|$ and
$|x(y_{1})-x(y_{2})|\leq 2|y_{1}-y_{2}|$ if
$y_{1},y_{2}\in\hat\Gamma$.
 In particular, if $|y-y(z)|\leq r $,
then $|x(y)-z|\leq 2r$,
so that $\chB^{+}_{r }(z)\subset B_{2r}(z)$ and
$$
\chB^{+}_{r }(z)\subset B_{2r}(z)\cap\Gamma
\subset  B_{2r}(z)\cap\Omega.
$$
 which proves the first
inclusion in \eqref{6.2.4}. 

Furthermore,
if $|x-z|\leq r/2\leq \rho_{0}(\Omega)$ and $x\in\Omega$, then, since
$z\in B_{4\rho_{0}(\Omega)}$, 
$$
x\in B_{5\rho_{0}(\Omega)}\cap\Omega\subset  \Gamma.
$$
Then
$|y(x)-y(z)|\leq  r$ and $y(x)\in\hat\Gamma$,
that is, $y(x)\in B^{+}_{ r}(y(z))$ so that
$x\in x(B^{+}_{ r}(y(z)))$, which yields the second
inclusion in \eqref{6.2.4}.

To prove (ii), it suffices to note that
$$
|\chB^{+}_{\nu r}(z)|=|B^{+}_{\nu r}(y(z)|
\leq N\nu^{d}|B^{+}_{ r}(y(z)|=N\nu^{d}
|\chB^{+}_{ r}(z)|.
$$
The lemma is proved.     \qed

\begin{corollary}
                                  \label{corollary 2.14.1}
If $z\in\chB^{+}_{2\rho_{0}(\Omega)}$ and
$r\leq (1/2) \rho_{0}(\Omega)$, then for
  any measurable function~$g$  
\begin{equation}
                                       \label{11.12.2}
\dashint_{\chB^{+}_{r }(z)}
|g(x)|\,dx\leq N(d)\dashint_{B_{2r}(z)\cap\Omega}
|g(x)|\,dx.
\end{equation}
 
\end{corollary}

Indeed, the domain of integration on the right
is wider than the one the left owing to \eqref{6.2.4},
and
$$
N(d)|\chB^{+}_{r }(z)|\geq|\chB^{+}_{4r}(z)|
\geq|B_{2r}(z)\cap\Omega|
$$
in light of \eqref{6.2.5} and \eqref{6.2.4}.

Next, set
$$
\chC^{+ }_{R}=[0,R^{2})\times\chB^{+}_{R} 
$$
and for $\rho+r\leq 4\rho_{0}(\Omega)$ and $z=(t,x)$
such that $x\in \chB^{+}_{\rho}$ and $t\in\bR$ define
$$
\chC^{+}_{r}(z)=[t,t+r^{2})\times\chB^{+}_{r}(x).
$$

\begin{lemma}
                                 \label{lemma 8,17.2}  
There exist 
  $\bar\gamma=\bar\gamma(d,\delta)   
\in(0,1)$ and $\alpha_{0}=\alpha_{0}(\delta,d)\in(0,1)$
such that for any $\gamma\in(0,\bar\gamma]$ and
  $\alpha\in(0,\alpha_{0})$, whenever

\(i\,\)   $r,\rho>0$,
 $\nu\geq 12$, $\rho+\nu r\leq 4\rho_{0}(\Omega) $,   
$z_{0}\in \chC^{+}_{\rho}$,

\(ii\,\)
$u\in W^{1,2}_{p}(\chC^{+}_{\rho+\nu r} )$   
 and
  $u(t,x)=0$ if $x\in\partial\Omega$, 

\noindent we have
$$
I_{r}(z_{0}):=\Big(\dashint_{\chC^{+}_{r}(z_{0})}
\dashint_{\chC^{+}_{r}(z_{0})}
|D^{2}u(z_{1})-D^{2}u(z_{2})|^{\gamma}\,dz_{1}dz_{2}\Big)^{1/\gamma}
$$
$$
\leq N\nu^{(d+2)/\gamma}\Big(\dashint_{\chC^{+}_{\nu r}(z_{0})}
(|\partial_{t}u +F(D^{2}u)|^{d+1}+|Du |^{d+1})\,dz\Big)^{1/(d+1)}
$$
\begin{equation}
                                         \label{8,17.6}
+N(\nu^{1+(d+2)/\gamma} r+\nu^{-\alpha})
\Big(\dashint_{\chC^{+}_{\nu r}(z_{0})}
|D^{2}u |^{d+1}\,dz\Big)^{1/(d+1)},
\end{equation}
where  the constants
$N$ depend  only on $d,\alpha$,  
$M_{2}(\Omega)$,  and $\delta$.

\end{lemma}

Proof.  
 By the change of variables formula
we see that $I_{r}(z_{0})$ equals
$$
 \Big(\dashint_{C^{+}_{r}(t_{0},y(x_{0}))}
\dashint_{C^{+}_{r} (t_{0},y(x_{0} ))}
|(D^{2} u)(x(z_{1}))-(D^{2}u)(x(z_{2}))|^{\gamma}
\,dz_{1}dz_{2}\Big)^{1/\gamma}.
$$

Then with
$A(y):= \partial x(y)/\partial y$ we define
$$
A=A(y(z_{0})),\quad
\check F(\sfu'')=F((A^{-1})^{*}\sfu''A^{-1}).
$$
As is easy to see, $D_{\sfu'' }\check F
\in\bS_{\check\delta}$, where $\check\delta=
\check\delta(d,\delta)\in(0,1]$.

Next,
introduce  the function
$$
\hat{u}(t,y)=u(t,x(y)), 
$$
 which belongs to $W^{1,2}_{p}(C^{+}_{\rho+\nu r})
  $, and, since $|y(x_{0})|< \rho $,
it also belongs to $W^{1,2}_{p}(C^{+}_{\nu r}
(t_{0},y(x_{0}))) $
and vanishes on $\partial_{x^{1}}C^{+}_{\nu r}
(t_{0},y(x_{0}))$
if this set in nonempty.
By Lemma \ref{lemma 8,17.3},  
since $\nu\geq 12$, we have
$$
\Big(\dashint_{C^{+}_{r} (t_{0},y(x_{0} ))}
\dashint_{C^{+}_{r}(t_{0},y(x_{0}))}
|D^{2} \hat{u}(z_{1})-D^{2} \hat{u}
(z_{2})|^{\gamma}\,dz_{1}dz_{2}\Big)^{1/\gamma}
$$
$$
\leq N\nu^{(d+2)/\gamma}
\Big(\dashint_{C^{+}_{\nu r}(t_{0},y(x_{0}))}
|\partial_{t}\hat u+
\check F(D^{2} \hat{u})|^{d+1}\,dz\Big)^{1/(d+1)}
$$
\begin{equation}
                                              \label{10.10.20}
+N\nu^{-\alpha}\Big(\dashint_{C^{+}_{\nu r}(t_{0},y(x_{0}))}
|D^{2} 
\hat{u}|^{d+1}\,dz
\Big)^{1/(d+1)}.
\end{equation}

Observe also  that for $y=y(x)$ and $x=x(y)$
$$
D \hat{u}(t,y)=(D u)(t,x)
A(y) ,
$$
where the $D $'s  are row vectors, and
\begin{equation}
                                              \label{6.2.01}
D^{2} \hat{u}(t,y)=A^{*}(y)
[D^{2} u(t,x)]A(y)
+[D_{k}u(t,x)]D^{2} x^{k}(y).
\end{equation}
Since   
$$
|A-A(y)|\leq N|y-y(x_{0})|,
$$
where $N$ depends only on $d$ and the bound on
$|D^{2}\psi|$, for $z_{i}=(t_{1},y_{i})
\in C^{+}_{r} (t_{0},y(x_{0} ))$, $i=1,2$, we have
$$
|D^{2} \hat{u}(z_{1})-D^{2} \hat{u}
(z_{2})|\geq(1/N)|D^{2} u(t_{1},x_{1})-D^{2} u
(t_{2},x_{2})|
$$
$$
-Nr (|D^{2} u(t_{1},x_{1})|+|D^{2} u
(t_{2},x_{2})|)-N(|D u(t_{1},x_{1})|+|D u(t_{2},x_{2})|),
$$
where $x_{i}=x(y_{i})$ and $N$ depends only on $M_{2}(\Omega)$
and $d$.  
Hence, the left-hand side
of \eqref{10.10.20} is greater than or equal to
$$
(1/N)I_{r}(z_{0})-N \Big(\dashint_{ \chC^{+}_{r}(z_{0})}
 (r|D^{2} u|+|D u|)^{\gamma}\,dz\Big)^{1/\gamma}   
$$
$$
\geq (1/N)I_{r}(z_{0})-Nr\Big(\dashint_{ \chC^{+}_{r}(z_{0})}
  |D^{2} u|^{d+1}\,dz\Big)^{1/(d+1)} 
$$
$$
-N\Big(\dashint_{ \chC^{+}_{r}(z_{0})}
 |D  u|^{d+1}\,dz\Big)^{1/(d+1)}   \geq (1/N)I_{r}(x_{0})
$$
\begin{equation}
                                    \label{6.2.3}
-N\nu r\Big(\dashint_{ \chC^{+}_{\nu r}(z_{0})}
  |D^{2} u|^{d+1}\,dz\Big)^{1/(d+1)} 
-N\nu\Big(\dashint_{ \chC^{+}_{\nu r}(z_{0})}
 |D  u|^{d+1}\,dz\Big)^{1/(d+1)}  ,
\end{equation}
where the first inequality follows by H\"older's
inequality and the second one is true owing to
\eqref{6.2.5}.  

In what concerns the first term on the right-hand side of
\eqref{10.10.20}, observe that, owing to the Lipschitz
continuity of $F$, the fact that
 $ |A(y)A^{-1}-(\delta^{ij})|\leq
N|y-y(x_{0})|$, 
and \eqref{6.2.01}, 
 we have (with $x=x(y)$)
$$
|\check F(D^{2} \hat u(y))-F(D^{2}u(x ))|
\leq
|F((A^{-1})^{*}A^{*}(y)
[D^{2} u(x)]A(y)A^{-1})
$$
$$
-F(D^{2}u(x ))|+
N|D u(x)|\leq N|y-y(z)||D^{2} u(x)|  
+N|D u(x)|.
$$
This and an easy estimate of the last term in 
\eqref{10.10.20} shows that its right-hand side is less
than
$$
N\nu^{(d+2)/\gamma}
\Big(\dashint_{\chC^{+}_{\nu r}(z_{0})}
|\partial_{t}u+ F(D^{2} u)|^{d+1}\,dx\Big)^{1/(d+1)}
$$
$$
+
N(\nu^{1+(d+2)/\gamma}r+\nu^{-\alpha})
\Big(\dashint_{\chC^{+}_{\nu r}(z_{0})}
|  D^{2} u |^{d}\,dx\Big)^{1/d}
$$
$$
+N(\nu^{(d+2)/\gamma}+\nu^{-\alpha})
\Big(\dashint_{\chC^{+}_{\nu r}(z_{0})}
|  D  u  |^{d}\,dz\Big)^{1/d}.
$$
 
Upon combining this result with what was said about
\eqref{6.2.3} we come to  
\eqref{8,17.6}. The lemma is proved.         \qed
 
Change of variables help derive
  Lemma \ref{lemma 8,17.2} from its ``flat''
counterpart. We also allude to it 
in the following remark.

\begin{remark}
                                    \label{remark 8,18.2}     
Suppose that Assumptions \ref{assumption 7,16.1} ,  
  and 
\ref{assumption 7,16.5}    and condition \eqref{1,19,500}
are satisfied.          Let $r\leq 4\rho_{0}(\Omega)$,  
  $p\geq d+1$, and $u\in W^{1,2}_{p}(\chC^{+}_{r})$
be such that   
  $u(t,x)=0$ if $x\in\partial\Omega$.          
Then  
$$
 \int_{\chC^{+}_{r} }
(|D^{2}u|^{\gamma}+|Du|^{\gamma})
\,dz \leq   
Nr^{d+2-\gamma (d+2)/p}\|\partial_{t}u+H[u]\|_{L_{p}
(\chC^{+}_{r } )}^{\gamma}
$$
$$
 +Nr^{d+2-\gamma (d+2)/p}\|\bar G\|_{L_{p}
(\chC^{+} _{r } )}^{\gamma} 
 +Nr^{d+2-2\gamma}\sup_{\partial' \chC^{+} _{r }  }|u|^{\gamma} ,
$$
where $N$ depend only on 
$\delta$, $K_{0}$, $d$, $p$,    and $M_{2}(\Omega)$
and the range of $\gamma$ is specified below.

Indeed, by using the notation from the above proof
and using equation \eqref{11.10.3} in
Lemma \ref{lemma 11.10.1} introduce the operators
$$
\cL u(t,x)=a^{ij}(t,x)D_{ij}u(t,x)+  b^{i}(t,x)D_{i}u(t,x),
\quad
\hat \cL \hat u(t,y)=[\cL u](t,x(y)).
$$
The operator $\hat \cL$ can be written as a differential
operator with derivatives with respect to $y$.
Clearly, its matrix of second-order derivatives
will belong to $\bS_{\hat \delta}$ for a $\hat \delta
=\hat \delta(\delta, M_{2}(\Omega))\in(0,1)$
and the drift term by magnitude will be dominated
by $N=N(K_{0},d,M_{2}(\Omega))$.
Since 
$$
|\partial_{t}\hat u(t,y)+\hat \cL \hat u(t,y)|\leq
 \big|\partial_{t}u(t,x(y))+H[u](t,x(y))\big|+\bar G(t,x(y))+K_{0}
|u(t,x(y))|
$$
 in $C^{+}_{4\rho_{0}(\Omega)}$,   
by Theorem \ref{theorem 8.16.1}
for an appropriate 
$\bar\gamma=\bar\gamma(d,\delta,K_{0}, M_{2}(\Omega))\in(0,1)$
and $\gamma\in(0,\bar\gamma ]$,
after using   scalings
and H\"older's inequality (to replace $d+1$ with $p$) we get,
$$
\int_{C_{r}^{+}}
(|D^{2}_{y}\hat u|^{\gamma}
+|D_{y}\hat u|^{\gamma})\,dydt
\leq Nr^{(d+2)(1-\gamma /p)}\Big(\int_{C_{r}^{+}}\big|
 \partial_{t}u+ H[u]\big|^{p}(t,x(y))\,dydt\Big)^{\gamma/p}
$$
$$
+Nr^{(d+2)(1-\gamma /p)}\Big(\int_{C_{r}^{+}}\big|
\bar G(t,x(y))\big|^{p}\,dydt\Big)^{\gamma/p}
+Nr^{d+2-2\gamma}\sup_{\partial' \chC^{+} _{r }  }|u|^{\gamma}.
$$
Now our assertion  follows after changing variables.
 
\end{remark}

\mysection
{A priori estimates in $W^{1,2}_{p}$
near the boundary and the proof of Theorem
\protect\ref{theorem 8,17.3}}
 
                                 \label{section 8,19.1}

We   assume that
$$
p>d+1,\quad
0\in\partial\Omega
$$
and
  take $ \rho_{0}=
\rho_{0}(\Omega)$, $\chC^{+}_{r} $, $\chC^{+}_{r}(z)$ from 
Section \ref{section 8.19.2}
and suppose that the assumptions of Theorem \ref{theorem 8,17.3}
are satisfied with  $\hat\theta$ and $\theta$ which are yet to be specified.
 
First we note the following.

\begin{lemma}
                                   \label{lemma 8,19.1}
 For any $q\in[1,\infty)$ and $\mu>0$ 
there exists  $\theta=\theta(d,\delta,K_{F},\mu,q)>0$
such  that,
if Assumption \ref{assumption 8,18.2} is satisfied
 with this $\theta$,
then
for any $\sfu'_{0}\in\bR$, $
z_{0}\in \chC^{+} _{2\rho_{0}(\Omega)}$ and
$2r\leq \rho_{0}(\Omega)\wedge R_{0}$, we have
$$
 \dashint_{\chC^{+}_{r}(z_{0})}\sup_{
\substack{\sfu''\in\bS,\\|\sfu''|>\tau_{0}}}
\frac{|F(\sfu'_{0},\sfu'',z)-\bar{F}(\sfu'')|^{q}}{| \sfu''|^{q}}\,dz\leq \mu^{q},
$$
where $\bar F=\bar F_{z,r, \sfu'_{0}}$ is taken from
Assumption \ref{assumption 8,18.2}.
\end{lemma}

For the proof of this lemma note that, in light
of Corollary \ref{corollary 2.14.1} and
Assumption \ref{assumption 8,18.2},  
for any $\sfu''\in\bS$ with $|\sfu''|=1$ we have  
$$
\dashint_{  \chC^{+}_{r }(z_{0})}\sup_{\tau>\tau_{0}}\tau^{-1}
|F(\sfu'_{0},\tau u'' , z)-\bar{F}(\tau u'')| \,dz\leq N(d)\theta
$$
if $2r\leq \rho_{0}(\Omega)\wedge R_{0}$. After that, as in the case
of Lemma \ref{lemma 8,15.1},
the assertion of the current lemma
 is obtained by repeating  
the proof of Lemma 5.1 of \cite{Kr_13}.  

Recall that $ \omega_{u}(\Pi,\rho )$
is introduced in Definition \ref{definition 1,18,1}.

\begin{lemma}
                                           \label{lemma 8,19.2}
Let $r,\rho\in(0,\infty)$ and
 $\nu\geq 12$ satisfy  $\rho+\nu r\leq
 4\rho_{0}(\Omega) $ and $\nu r\leq
 R_{0}  $.
Take   
$$
\mu\in(0,\infty), \quad\beta\in(1,\infty),
$$ 
and suppose that Assumption \ref{assumption 8,18.2}
is satisfied with $\theta=\theta(d,\delta,K_{F},\mu,\beta d+\beta)$
\(see Lemma \ref{lemma 8,19.1}\).
Assume that   we are given a function  
  $u\in W^{1,2}_{p}(\chC^{+} _{\rho+\nu r})$
and $u(t,x)=0$ if $x\in
\partial \Omega$.
Use $I_{r}(z_{0})$ introduced in \eqref{8,17.6}.
 
 Then,  for  $\gamma$ and $\alpha$ from Lemma \ref{lemma 8,17.2}, for  
$z_{0}\in\chC^{+} _{  \rho }$,
we have
$$
I_{r} (z_{0})\leq 
N\eta
\big(
\dashint_{\chC^{+}_{\nu r}(z_{0})}
|D^{2}u |^{\beta' (d+1)}\,dz\Big)^{1/(\beta'(d+1))}  
$$
$$
+N\nu^{(d+2)/\gamma}\Big(
\dashint_{ \chC^{+}_{\nu r}(z_{0})}(\big|\partial_{t}u+F[u]\big|^{d+1}
+|Du |^{d+1})\,dz\Big)^{1/(d+1)}+N\tau_{0}\nu^{(d+2)/\gamma},
$$
where
$$
\eta=\big(\mu+\nu  r
+\omega_{F,u,\chC^{+} _{\rho+\nu r}}(\nu r )\big)\nu^{(d+2)/\gamma}+\nu^{-\alpha},
$$
and the constants $N$ depend only 
on $  d,p, 
  K_{F},\delta$, and   $M_{2}(\Omega)$.
\end{lemma}

The proof of this lemma is based on Lemma
\ref{lemma 8,19.1} and, in light of
Lemma \ref{lemma 8,17.2}, is practically identical to that 
of Lemma \ref{lemma 7,29.3}.

We now come to the main   a priori estimate
 near the boundary for 
nonlinear parabolic equations with VMO ``coefficients".

\begin{theorem}
                               \label{theorem 8,19.2}
Take $p>d+1$, 
let $R >0$ satisfy
$$
2R\leq \rho_{0}(\Omega)\wedge R_{0},  
$$
and let $u\in W^{1,2}_{p}(\chC^{+} _{2R})$ be
such that $u(t,x)=0$ if $ x\in 
\partial \Omega$.
Then there exist  constants $\hat\theta,\theta\in(0,1]$, depending
only on $d,p,\delta$, and $K_{F}$, such that
if Assumptions \ref{assumption 7,16.01}
and \ref{assumption 8,18.2}   
 are satisfied with these $\hat\theta$ and $\theta$,
respectively,
then there exist   constants $N$, depending only
on $R_{0},d,p, K_{0},K_{F},\delta$,   $\rho_{0}(\Omega)$,
  $M_{2}(\Omega)$, and
the function $\omega_{F,u,\chC^{+}_{2R}}$
\(see Definition \ref{definition 1,18,1}\),  
such that
$$
\|D^{2}u\|_{L_{p}(\chC^{+}_{R  })}\leq
N \|\partial_{t}u+H[u]\|_{L_{p}(\chC^{+}_{2R})} +N\| \bar G\|_{  
L_{p}(\chC^{+}_{2R})}
+N\tau_{0}  
$$
\begin{equation}
                                          \label{8,19.4}
 +N\| u\|_{L_{p}(\chC^{+}_{2R  })}
+
N  R  ^{-\chi }
  \|
\,|D^{2}u|^{\gamma}
\|^{1/\gamma}_{L_{1}(\chC^{+}_{2R })},  
\end{equation}

$$
\|D^{2}u\|_{L_{p}(\chC^{+} _{R })}
\leq  N \|\partial_{t}u+H[u]\|_{L_{p}(\chC^{+} _{2R})}
+N \|\bar G\|_{L_{p}(\chC^{+} _{2R})}
$$
\begin{equation}
                                           \label{8,19.5}
+N\tau_{0}
+N R   ^{(d+2)/p-2}\sup_{\chC^{+} _{2R}}|u|,
\end{equation}
where    $\chi =
 (d+2)(1/\gamma-1/p )$
and $\gamma$ is the same as in Lemma \ref{lemma 8,19.2}.
 
\end{theorem}

Proof. Whenever it makes sense, for
  $\rho\leq\rho_{0}(\Omega)$,
and $z\in\chC^{+ }_{2\rho_{0}(\Omega)}$ introduce
$$
 h^{\vsharp}_{\Omega,\gamma,\rho}(z)=
\sup\{I_{r }(h,z_{0}):z_{0}\in\bR^{+}\times\Omega,
r\in(0,\rho],\chC^{+}_{r}(z_{0})\ni z \},
$$
where
$$
I_{r}(h,z_{0})=\bigg(
\dashint_{ \chC^{+}_{r}(z_{0})}\dashint_{\chC^{+}_{r}(z_{0})}
|h(z_{1})-h(z_{2})|^{\gamma}\,dz_{1}dz_{2}\bigg)^{1/\gamma}.
$$
The reader should pay attention to the above curved
sharp symbol, reminding of curved boundaries.

Observe that, if $r\leq \rho\leq\rho_{0}(\Omega)$ and $z\in
\chC^{+ }_{2\rho_{0}(\Omega)}
\cap \chC^{+}_{r}(z_{0})$, then   $\chC^{+}_{r}(z_{0})
\subset \chC^{+ }_{4\rho_{0}(\Omega)}$, so that 
$h^{\vsharp}_{\Omega,\gamma,\rho}(z)$ is well defined on 
$\chC^{+ }_{2\rho_{0}(\Omega)}$
even if $h$ is given only on $\chC^{+ }_{4\rho_{0}(\Omega)}$
($\subset\Omega$).
 
Then take $\varepsilon\in(0,1]$ to be specified later,
 take $R_{1}<R_{2}\leq 2R$ such that 
$$
 R_{2}\leq 2R_{1},\quad
  R_{2}-R_{1}\leq \varepsilon  R_{0},
$$
take $\nu\geq 12$, and set
$$
r_{0}=(R_{2}-R_{1})/(\nu+1),\quad
\kappa=r_{0}/R_{1}\quad(\leq(R_{2}-R_{1})/(2R_{1})\leq1/2).
$$

We are going to use Theorem \ref{theorem 8,12.1}
according to which, if  $h\in L_{p}( \chC^{+}_{R_{2}})$,
then
\begin{equation}
                                             \label{4,26,1}
\|h\|_{L_{p}( \chC^{+}_{R_{1}})}
\leq N\|h^{\vsharp}_{\Omega,\gamma,r_{0}}\|
_{L_{p}( \chC^{+}_{R_{1}})}+ N
\nu^{\chi_{1}}(R_{2}-R_{1})^{-\chi_{1}}
R_{1}^{\chi_{1}-\chi_{2}}
\|\,|h|^{\gamma}\|^{1/\gamma}_{L_{1}( \chC^{+}_{R_{1}})},
\end{equation}
where $\chi_{1}=(d+2)/\gamma$, $\chi_{2}=(d+2)(1/\gamma-1/p)$,
and the constants
 $N$ depend only on $d,\gamma$, and $p$.

Next for $z\in \chC^{+}_{2R}$ ($\subset \chC^{+}_{\rho_{0}(\Omega)}$)
 define
$$
\bM_{\Omega} h(z) =
\sup\Big\{\dashint_{ \chC^{+}_{r}(z_{0} )}
|h(y)|\,dy:2r\leq \rho_{0}(\Omega),z_{0}\in\bR^{+}\times\Omega,
\chC^{+}_{r}(z_{0} )\ni z\Big\}.
$$
 
Observe that, owing to the fact that $\Omega\in C^{1,1}$   
and to Corollary \ref{corollary 2.14.1}, 
if $z_{0}\in \chC^{+}_{2\rho_{0}(\Omega)}$ 
and
$r\leq (1/2) \rho_{0}(\Omega)$,
$$
\dashint_{ \chC^{+}_{r}(z_{0})}
|h |\,dy\leq   N\dashint_{C_{2r}(z_{0}) }
|h |I_{\chC^{+}_{ 2R } }\,dy,
$$
 where  
 $N$ depends only on $d$, $\rho_{0}(\Omega)$, and $M_{2}(\Omega)$. Therefore,
for $z\in\chC^{+}_{2R}$
\begin{equation}
                                            \label{8,20.1}
\bM_{\Omega} h(z)\leq N\bM  
hI_{\chC^{+}_{2R} }(z).
\end{equation}
 
The above conclusion \eqref{8,20.1} is, actually,
 also based on the fact similar to the following.
For $r\leq r_{0}$, $z\in\chC^{+}_{R_{1}}$,
and $z_{0}\in\bR^{+}\times\Omega$, such that 
  $\chC^{+}_{r}(z_{0})\ni z$,
we have $z_{0}\in \chC^{+}_{\rho}$, where $\rho=R_{1}+r $.
In this situation also $\rho+\nu r\leq R_{2}<4\rho_{0}(\Omega)$ and
$\nu r\leq\varepsilon  R_{0}$ and it follows from 
 Lemma \ref{lemma 8,19.2} that
$$ 
I_{r}(z_{0}) \leq
N\nu^{(d+2)/\gamma}\bM^{1/(d+1)}_{\Omega}\big(\big|
\partial_{t}u+F[u]\big|^{d+1}I_{\chC^{+}_{R_{2}}}
 \big)(z)+N\tau_{0}\nu^{(d+2)/\gamma}  
$$
$$
+N\eta
\bM_{\Omega}^{1/(\beta'(d+1))}
\big(|D^{2}u|^{\beta'(d+1)}I_{  \chC^{+}_{R_{2}}}\big)
(z)
$$
$$
+N\nu^{(d+2)/\gamma}\bM_{\Omega}^{1/(d+1)}
\big(|Du|^{d+1}I_{\chC^{+}_{R_{2}}}\big)(z),
$$
where 
$$
\eta=\big(\mu+\nu  r_{0}
+\omega_{F,u,\chC^{+} _{R_{2}}}(\varepsilon  R_{0})\big)\nu^{(d+2)/\gamma}+\nu^{-\alpha},
$$  
By definition and \eqref{8,20.1}
we obtain that on $ \chC^{+}_{R_{1}}$
$$
(D^{2}u)^{\vsharp}_{\Omega,\gamma,r_{0}}
\leq N\nu^{(d+2)/\gamma}\bM ^{1/(d+1)}
\big(\big|\partial_{t}u+F[u]\big|^{d+1}I_{\chC^{+}_{R_{2}}
 }\big)+N\tau_{0}\nu^{(d+2)/\gamma}
$$
$$
+N\eta
\bM ^{1/(\beta'(d+1))}\big(|D^{2}u|^{\beta'(d+1)}I_{ \chC^{+}_{R_{2}}}\big)
$$
$$
+N\nu^{(d+2)/\gamma}\bM ^{1/(d+1)}\big(|Du|^{d+1}I_{\chC^{+}_{R_{2}}
 }\big).
$$

Thanks to \eqref{4,26,1} and
the Hardy-Littlewood maximal function theorem,
  by taking $\beta$ so that $p
>\beta'd$, we obtain
$$
\|D^{2}u\|_{L_{p}(\chC^{+}_{R_{1}})}
\leq N \nu^{(d+2)/\gamma}\big\|\partial_{t}u+
F[u]\big\|_{L_{p}(\chC^{+}_{R_{2}})}+
N\tau_{0}\nu^{(d+2)/\gamma}
$$
$$
+\Big[N \big(\mu+\nu  r_{0}
+\omega_{F,u,\chC^{+} _{R_{2}}}(\varepsilon  R_{0})\big)\nu^{(d+2)/\gamma}+N_{0}\nu^{-\alpha}
\Big]
\|D^{2}u\|_{L_{p}(\chC^{+}_{R_{2}})}
$$
$$
+N\nu^{(d+2)/\gamma}\|D u\|_{L_{p}(\chC^{+}_{R_{2}})}
+N\nu^{\chi_{1}}(R_{2}-R_{1})^{-\chi_{1}}
R_{1}^{\chi_{1}-\chi_{2}}
\big\|\,|D^{2}u\big|^{\gamma}\|^{1/\gamma}_{L_{1}(\chC^{+}_{2R  })},
$$
where the constants $N$, $N_{0}$ depend only on $ d$, $p$, $K_{F}$, 
and $\delta$.

This estimate looks almost like \eqref{4,21,3}.
Then we repeat the argument after \eqref{4,21,3}
and choose and fix $\varepsilon$, $\nu$, $\hat\theta$, and $\mu$, recall   
what $r_{0}$ is, and conclude that
$$
\|D^{2}u\|_{L_{p}(\chC^{+}_{R_{1}})}
\leq N \big\|\partial_{t}u+
H[u]\big\|_{L_{p}(\chC^{+}_{R_{2}})}+
N\tau_{0}\nu^{d/\gamma}
$$
$$
+ (5/8+N (R_{2}-R_{1}) )
\|D^{2}u\|_{L_{p}(\chC^{+}_{R_{2}})}
+N\|u\|_{L_{p}(\chC^{+}_{2R})}+N\|\bar G\|_{L_{p}(\chC^{+}_{2R  })}
$$
$$
+N(R_{2}-R_{1})^{-\chi_{1}}
R_{1}^{\chi_{1}-\chi_{2}}
\big\|\,|D^{2}u\big|^{\gamma}\|^{1/\gamma}_{L_{1}(\chC^{+}_{2R  })}.
$$
 
After that, to prove \eqref{8,19.4}, it suffices to repeat 
almost literally what follows \eqref{8,20.3} (only replacing $C$
with $\chC$).
By using Remark \ref{remark 8,18.2}   we estimate the   
last term in \eqref{8,19.4} and then finish the proof of the theorem
in the same way as in the case of Lemma \ref{lemma 7,30.1}.  
The theorem is proved.    \qed

{\bf Proof of Theorem \ref{theorem 8,17.3}}.   
To start, assume that $g\equiv0$.
Observe that in that case
 we may assume that $u(t,x)$ is defined for
$t\geq T$, $x\in\bar\Omega$, as zero and still satisfies there
\eqref{8,17.7}. It suffices for the
 latter  that $H(0,t,x)=0$ if $t\geq T$, which is easy
to accommodate without altering our assumptions
just by replacing $G(\sfu,t,x)$ and $\bar G$ with $G(\sfu,t,x)I_{t<T}$
and $\bar G I_{t<T}$, respectively.
After such extension $u\in W^{1,2}_{p}(\bR_{+}\times\Omega)$.

Take $\hat\theta$ and $\theta$ which suit  both Lemma   
\ref{lemma 7,30.1}
and  Theorem \ref{theorem 8,19.2}, and take  
 $$
 R=(1/2)(\rho_{0}(\Omega)\wedge R_{0}).
$$

Theorem \ref{theorem 8,19.2} allows us to
estimate the $W^{1,2}_{p}$-norm of $u$ in the domain
$\chC^{+}_{R}=\chC^{+}_{R} $ associated
with the origin, that is assumed to belong to
$\partial\Omega$. Of course, one can take any point  
$z_{0}=(t_{0},x_{0})\in[0,\infty)\times\partial\Omega$ as the origin and
apply Theorem \ref{theorem 8,19.2} to
$\chC^{+}_{R}(z_{0})$ and $\chC^{+}_{2R}(z_{0})$ in place of
$\chC^{+}_{R}$ and $\chC^{+}_{2R}$, respectively,
where by $\chC^{+}_{\rho}(z_{0})$ we, naturally, mean
the sets 
$$
[t_{0},t_{0}+\rho^{2})\times\chB^{+}_{\rho}(x_{0}),
$$ 
where $\chB^{+}_{\rho}(x_{0})$ is constructed in Lemma
\ref{lemma 6.2.2} but with $x_{0}$ in place of $0$
and relative to the coordinate system $\Psi(x_{0})$
associated with $x_{0}$ as described before
Lemma \ref{lemma 2.14.1}.
According to that, we  find
  finitely many 
$$
z_{i}\in [0,T+ R^{2}]\times\partial\Omega
$$
and $\rho>0$ depending only on $\diam(\Omega)$, $\rho_{0}(\Omega)$,
$M_{2}(\Omega)$, and $T$ such that
$$
 \bigcup_{i}\chC^{+}_{R}(z_{i})\cup\Big([0,S-\rho^{2})\times    
\Omega^{\rho}\Big)\supset\Pi,
$$
where $S=T+R^{2}$.

By Theorem \ref{theorem 8,19.2}, for any $i$
(recall that $\bar G(t,x)=0$ for $t\geq T$)
$$
\|D^{2}u\|^{p}_{L_{p}(\chC^{+} _{R }(z_{i}))}
\leq  N  \|\bar G\|^{p}_{L_{p}(\Pi)}
+N\tau^{p}_{0}
+N \sup_{\Pi }|u|.
$$
By Lemma \ref{lemma 7,30.1}
$$
\|D^{2}u\|^{p}_{L_{p}(0,S-\rho^{2})\times
\Omega^{\rho})}
\leq  N  \|\bar G\|^{p}_{L_{p}(\Pi )}
+N\tau^{p}_{0}
+N \sup_{\Pi}|u|.
$$
We sum up these estimates and come to
\eqref{8,18.1}.
This proves the theorem
if $g\equiv0$.

In the general case introduce 
$\hat g(z)=(g(z),Dg(z),D^{2}g(z))I_{\Pi}(z)$
and
$$
\hat H(\sfv,z)=H(\sfv+\hat g(z),z)
  +\partial_{t}g(z)I_{\Pi}(z)  ,\quad w(z)=u(z)-g(z).
$$
Observe that $\hat H[w]=0$ in $\Pi$ (a.e.) and
$w\in W^{1,2}_{p}(\Pi)$ and $w=0$ on $\partial'\Pi$.
 Furthermore, for
$$
 \hat F(\sfv'_{0}, \sfv'',z)=F(\sfv'_{0}+
\hat g'_{0} ,\sfv'',z) ,\quad
\hat G(\sfv,z):=\hat H(\sfv,z)-
 \hat F(\sfv'_{0}, \sfv'',z) 
$$
  one easily obtains that $|\hat G(\sfv,z)|
\leq\hat\theta|\sfv''|+K_{0} |\sfv'| +\bar{\!\hat G}$, where
$$
\bar{\!\hat G}=|\partial_{t}g|I_{\Pi}+ N|D^{2}g |I_{\Pi}+\bar G+K_{0}
 \big( g  ^{2}+|Dg|^{2} \big)^{1/2}I_{\Pi}
 $$
with $N$ depending only  on $K_{F}$ and $d$.

Also for $\sfv'_{0}\in\bR$,  $r\in(0, R_{0}]$, and $z\in \Pi$   
we set
$$
\bar{\!\hat F} (\sfv'')=
\bar{\!\hat F}_{z,r, \sfv'_{0}}(\sfv'')=\bar F_{z,r, \sfv'_{0}+
\hat g'_{0}(z)}(\sfv''),
$$
take $\theta_{0}=\theta_{0}\big(d,p,\delta, K_{F},M_{2}(\Omega)\big)$
defined above in the first part of the proof where $g\equiv0$,
 find $\hat R_{0}\leq R_{0}$ such that $\omega_{F,g,\Pi}(\hat R_{0})
\leq \theta_{0}/2$ and then require the original
Assumption \ref{assumption 8,18.2} (iii) to be satisfied
with $\hat R_{0}$ and $\theta_{0}/2$ in place of $R_{0}$  
and $\theta$, respectively.

Then we see   that the above result is applicable to $w$,
and along with the embedding inequality:
$|g|\leq N\|g\|_{W^{1,2}_{p}(\Pi)}$,
  lead to \eqref{8,18.1} in the general case.
The theorem is proved.               \qed

\mysection
{Proof of 
Theorem \protect\ref{theorem 8,18.2}}
                                       \label{section 8,26.1}

The proof of Theorem \ref{theorem 8,18.2}  is based on the
following.

\begin{theorem}
                                        \label{theorem 4,5,1}
  Suppose that Assumption \ref{assumption 7,16.5}
is satisfied,  the number
$$
\bar{H} :=\sup_{\sfu', t, x}\big(|H (\sfu',0,t, x)|-K_{0}|\sfu'|\big)
\quad(\geq0)
$$
is finite, and
$g\in W^{1,2}_{\infty }(\bR^{d+1})$.

Then there exists a convex positive homogeneous of degree
one function $P(\sfu'')$ such that 
  at all points of its differentiability 
$D_{\sfu''}P \in
\bS_{\bar\delta}$, where $\bar\delta=\bar\delta(d,\delta)\in(0,\delta)$,
and for $P[u]=P(D^{2}u)$
and any $K>0$ the equation
\begin{equation}
                                               \label{9.28.4}
\partial_{t}v+ \max(H[v],P[v]-K)=0  
\end{equation}
in $\Pi$ with boundary condition $v=g$ on $\partial'\Pi$
has a solution  $v\in
 W^{1,2}_{p}(\Pi)$ for any $p\geq 1$.
\end{theorem}

This theorem follows from Theorem 2.1 of \cite{Kr_13.1},
proved there under the additional conditions that $\Omega\in C^{2}$
and that there is an increasing
 continuous function $\omega(r)$, $r\geq0$,
such that $\omega(0)=0$ and
$$
|H(\sfu',\sfu'',t,x)-H(\sfv',\sfu'',t,x)|\leq\omega(|\sfu'-\sfv'|)
$$
for all $\sfu,\sfv,t$, and $x$. That these additional conditions can be dropped
will be proved elsewhere.

{\em Step 1\/}. 
We take $P(\sfu'')$   from Theorem \ref{theorem 4,5,1},  
and
first we assume that $g\in W^{1,2}_{\infty}
(\Pi)$
and there exists   constants $N_{0},\bar H$ such that,
for all $t,x,\sfu'$,
\begin{equation}
                                                    \label{1,22,2}
|H(\sfu',0,t,x)|\leq N_{0}|\sfu'|+\bar H.
\end{equation}

By Theorem \ref{theorem 4,5,1} 
for any $K>0$ there 
 exists a function $v_{K}$ which is  in $W^{1,2}_{p}(\Pi)$  
for any $p>1$, such that 
  $v_{K}=g$
on $\partial'\Pi $, and it satisfies 
\begin{equation}
                                              \label{1,22,3}
\partial v_{K}+H_{K}[v_{K}]=0
\quad\text{in}\quad \Pi \,\,\text{\rm (a.e.)},
\end{equation}
where
$$
H_{K}(\sfu,t,x)=\max(H(\sfu,t,x) ,P(\sfu'')-K).
$$
 
Set 
$$
F_{K}(\sfu'_{0},\sfu'',t,x)
=\max(F(\sfu'_{0},\sfu'',t,x),P( \sfu'')-K),
$$
$$
\bar F_{K,z,r,\sfu'_{0}}(\sfu'')=
\max(\bar F_{z,r,\sfu'_{0}}(\sfu''),P(\sfu'')-K),
$$
$$
 G_{K}(\sfu,t,x)=H_{K}(\sfu,x)-F_{K}(\sfu'_{0},\sfu'',t,x).
$$

It is not hard to see that
 Assumptions
\ref{assumption 8,18.2}, \ref{assumption 7,16.01},  
and  
 \ref{assumption 7,16.5} are satisfied for $H_{K},F_{K}$,
and $G_{K}$ in place of $H ,F $,
and $G $, respectively, with the same $K_{0}$, $\bar G$, $R_{0}$,
$\theta$, $\hat\theta$, $\omega_{F}$, with $\bar\delta$ in place of $\delta$
and $\bar\delta^{-1}+K_{F}$ in place of $K_{F}$.
By Theorem \ref{theorem 8,17.3}
there exist  constants $\hat\theta,\theta\in (0,1]$, depending
only on $d$, $p$, $\delta$, $K_{F}$,  $ \rho_{0}(\Omega)$,
and $M_{2}( \Omega)$,
such that, if Assumptions \ref{assumption 8,18.2}  
and \ref{assumption 7,16.01} 
 are satisfied with these $\theta$ and
 $\hat\theta$, respectively,  then for any $K>0$, we have
$$
\|v_{K}\|_{W^{1,2}_{p}(\Pi)} \le N \big(
\| \bar G\|_{L_p(\Pi)} +\|g \|_{W^{1,2}_{p}(\Pi)} +
 \|v_{K}\|_{C(\Pi)}\big)
 +N \tau_{0} , 
$$
where the constants $N$  depend only on
 $K_{0}$, $K_{F}$, $d$, $p$, $\delta$,   
$R_{0}$, 
 $\diam(\Omega)$,
$\rho_{0}(\Omega)$, $M_{2}(\Omega)$, and the function 
$\omega_{F,v_{K},\Pi}$
(independent of $N_{0}$ and $\bar H$).

Since $H_{K}$ satisfies \eqref{1,19,500}, formula
\eqref{11.10.3} is valid with $v_{K}$ and $H_{K}$
in place of $u$ and $H$. This converts equation
\eqref{1,22,3} into a linear equation and by the
well-known results from the linear theory
allows us to estimate $|v_{K}|$ and the modulus of continuity
of $v_{K}$ through that of $g$, $\sup|g|$, and $\|\bar G\|_{L_{d+1}(\Pi)}$
with constants independent of $K$.

  Thus,
\begin{equation}
                                                \label{8,26.4}
\|v_{K}\|_{W^{1,2}_{p}(\Pi)}
\leq N \big(\|\bar{G}\|_{L_{p}(\Pi)}+\|g \|_{W^{1,2}_{p}(\Pi)}
\big)
  +N \tau_{0} ,
\end{equation} 
where the constants $N$ are independent of $K$. 
 
In this way we completed a crucial step consisting of
 obtaining a uniform
control of the $W^{1,2}_{p}(\Pi)$-norms of $v_{K}$.

Next, we  let $K\to\infty$. Estimate \eqref{8,26.4} guarantees that
there is a sequence $K_{n}\to\infty$
as $n\to\infty$ and $v\in W^{1,2}_{p }(\Pi)$
  such that $v_{K_{n}}\to v$ weakly in $W^{1,2}_{p}(\Pi )$
 and $v_{K_{n}}\to v$
  uniformly in $\bar \Pi$. Then, of course, $v=g$
on $\partial'\Pi$.
The said weak convergence implies  pointwise 
 convergence $Dv_{n}\to Dv$ in $\Pi$
in light of the compactness of the embedding $W^{1,2}_{p}\subset
C^{0,1}$ ($p>d+2$). 

Next, for $m=1,2,...$ define
$$
H^{m}(\sfu'',t,x)=\sup_{n\geq m}\max(H(v_{K_{n}}(t,x),
Dv_{K_{n}}(t,x),\sfu'',t,x), P(\sfu'')-K_{n}).
$$
Observe that $H^{m}(\sfu'',t,x)$ are Lipschitz continuous
in $\sfu''$ and at all points of differentiability satisfy
$D_{\sfu''}H^{m}\in\bS_{\bar\delta}$. 
Also
$$
|H^{m}(0,t,x)|\leq K_{0}\max_{n\geq m}\Big(|v_{K_{n}}(t,x)|+
|Dv_{K_{n}}(t,x)|\Big)+\bar G(t,x),
$$
which is in $L_{p,\loc}(\Pi)$.
Therefore, the operators $H^{m}[u]$
  fit into  the scheme of Section
3.5 of \cite{Kr_85}. Furthermore, for $n\geq m$ obviously
$$
\partial_{t}v_{K_{n}}+
H^{m}(v_{K_{n}},t,x)\geq0
$$
(a.e.) in $\Pi$. By Theorem  
 3.5.9  of \cite{Kr_85} we conclude that for any $m$
\begin{equation}
                                                \label{1,19,10}
\partial_{t}v+\sup_{n\geq m}\max(H(v_{K_{n}} ,
Dv_{K_{n}} ,D^{2}v ,t,x), P(D^{2} v)-K_{n})\geq0
\end{equation}
(a.e.) in $\Pi$. We fix $(t,x)$ at which \eqref{1,19,10}
holds for all $m$ (that is, we fix almost any $(t,x)$) and   since
$H(\sfu',\sfu'',t,x)$ is continuous in $\sfu'$, we have that
$$
|H(v_{K_{n}}(t,x),
Dv_{K_{n}}(t,x),D^{2}v(t,x),t,x)
$$
$$
-
H(v (t,x),
Dv (t,x),D^{2}v(t,x),t,x)|\to 0
$$
as $n\to\infty$. Then, in light of \eqref{1,19,10},
$$
\partial_{t}v(t,x)+
\max(H(v (t,x),
Dv (t,x),D^{2}v(t,x),t,x), P(D^{2}v(t,x))-K_{m})\geq o(1),
$$
which for $m\to\infty$ yields
$$
\partial_{t}v(t,x)+ H(v (t,x),
Dv (t,x),D^{2}v(t,x),t,x)=\partial_{t}v(t,x)+H[v](t,x)\geq0.
$$

The inequality $\partial_{t}v+H[v]\leq0$ is proved similarly
starting from the function
$$
\inf_{n\geq m}\max(H(v_{K_{n}}(t,x),
Dv_{K_{n}}(t,x),\sfu'',t,x), P(\sfu'')-K_{n}).
$$
  Owing to \eqref{8,26.4},
of course, $v\in W^{1,2}_{p}(\Pi)$ and 
\eqref{8,26.4} holds with $v$ in place of $v_{K}$.

This proves the theorem if condition
\eqref{1,22,2} is satisfied   and  $g\in W^{1,2}_{\infty}
(\Pi)$.
 {\em Step 2\/}. Assume that $g\in W^{1,2}_{\infty }(\bR^{d+1})$
and abandon \eqref{1,22,2}. 
Let 
$\eta(t)=t$ for $|t|\leq 1$ and $\eta(t)=\sign t$ for $|t|\geq1$.
For $n=1,2,...$ define $\eta_{n}(t)=n\eta(t/n)$ and
$$
\hat H^{n}(\sfu,t,x)=H(\sfu,t,x)-H(\sfu',0,t,x)+\eta_{n}(H(\sfu',0,t,x)),
$$
$$
\hat G^{n}(\sfu,t,x)=\hat H^{n}(\sfu,t,x)-F(\sfu'_{0},\sfu'',t,x).
$$

Then
$$
|\hat G^{n}(\sfu,t,x)|=| G (\sfu,t,x)+\eta_{n}(H(\sfu',0,t,x))-
H(\sfu',0,t,x)|
$$
$$
\leq \hat\theta|\sfu''|+2K_{0}|\sfu'|+2\bar G(t,x),
$$
so that Assumption \ref{assumption 7,16.01}
is satisfied for $\hat H^{n}$ with $2K_{0}$ and $2\bar G$ in place of
$K_{0}$ and $\bar G$. Assumptions \ref{assumption 8,18.2}
and \ref{assumption 7,16.5} are also valid for $\hat H^{n}$ 
with the same parameters.

Furthermore
$$
|\hat H^{n}(\sfu',0,t,x)|=|\eta_{n}(H(\sfu',0,t,x))|,
$$
which is bounded.

Hence there are $\hat\theta$ and $\theta$ as in Step 1,
 for any $n$, there exists
$u^{n}\in W^{1, 2}_{p }(\Pi )\cap C(\bar\Pi )$
 satisfying 
$$
\partial_{t}u^{n}+\hat H^{n}[u^{n}]=0
$$
 in $\Pi $ \(a.e.\)  and such
that $u=g $ on $\partial'\Pi $.
Estimate \eqref{8,18.1}, applicable to $v^{n}$ by the above
again guarantees that the $W^{1,2}_{p}(\Pi )$-norms
of $v^{n}$ are bounded  and $v^{n}$
are equicontinuous in $\bar\Pi$. This enables us to find
a subsequence
  $v^{n'}$ and a function $v\in W^{1,2}_{p }(\Pi)$ such that
$v^{n'}\to v$ weakly in $W^{1,2}_{p}(\Pi )$
  and $v^{n'}\to v$
  uniformly in $\bar \Pi$. Then, of course, $v=g$
on $\partial'\Pi$.

After that we repeat the rest of Step 1 by taking
$$
\sup_{n'\geq m}
\Big[H(v^{n'},Dv^{n'},\sfu'',t,x)-H(v^{n'},Dv^{n'},0,t,x)
+\eta_{n'}(H(v^{n'},Dv^{n'},0,t,x))\big]
$$
in place of $H^{m}(\sfu'',t,x)$. One thing which makes
the argument here easier is that for any $(t,x)\in\Pi$
$$
 - H(v^{n },Dv^{n },0,t,x)
+\eta_{n }(H(v^{n },Dv^{n },0,t,x))=0
$$
if $n$ is large enough.

In this way we finish Step 2. 
Finally, to treat the general $g\in W^{1,2}_{p}(\Pi)$
it suffices to use approximations and very simple
arguments about passing to the limit, which we
have seen already above. This step is left to
the reader. 
 The theorem is proved.    \qed

 \mysection{Appendix}
                                               \label{section 4,3,3}

Fix $\gamma\in(0,1]$ and for $r\in(0,\infty)$
and $z\in\bR^{d+1}$ define
\begin{equation}
                                                 \label{8.13.6}
I_{r}(h,z)=
\bigg(
\dashint_{C_{r}(z )}\dashint_{C_{r}(z  )}
|h(z_{1})-h(z_{2})|^{\gamma}\,dz_{1}dz_{2}\bigg)^{1/\gamma}
\end{equation}
whenever the right-hand side makes sense.

For   $\rho>0$
introduce the restricted sharp function of $h$   by
the formula
\begin{equation}
                                              \label{8,13.1}
 h^{\shharp}_{Q,\gamma,\rho}(z)=\sup\big\{I_{r}(  h,z_{0}):
z_{0}\in
Q,r\in(0,\rho],C_{r}(z_{0})\ni z\big\}
\end{equation}
whenever it makes sense.
Note that, if $Q=\bR_{+}\times\bR^{d}$,
 $h^{\shharp}_{Q,\gamma,\rho}$ is well defined
in $C_{R}$
for measurable $h$ even defined only in $C_{R+2\rho}$.
 
\begin{theorem}
                                     \label{theorem 8,4.1}
Let $p\in(1,\infty)$, $\kappa\in(0,1]$, $R\in(0,\infty)$,
 and $h\in L_{p}(C _{R(1+ 2\kappa)})$.
Let $Q=\bR_{+}\times\bR^{d}$.
Then
\begin{equation}
                                            \label{8,4.1}
\|h\|_{L_{p}(C _{R})}\leq N\big\|h^{\shharp}_{Q,\gamma,\kappa R}\big\|
_{L_{p}(C _{R })}+ N\kappa^{-\chi_{1}}
R^{-\chi_{2}}\big\|\,|h|^{\gamma}\big\|
_{L_{1}(C _{ R })}^{1/\gamma},
\end{equation}
where  $\chi_{1}=(d+2)/\gamma$, $\chi_{2}=
(d+2)(1/\gamma-1/p )$
 and the constants
 $N$ depend only on $d,\gamma$, and $p$.
\end{theorem}

This theorem will be proved elsewhere by closely following
the proof of Theorem 7.1 of \cite{Kr_13}
given there in the elliptic framework.

The remaining results of this section treat smooth  cylinders
or smooth domains.
If $\Omega\in C^{1,1}$ and $0\in\partial\Omega$,
 we assume  that the original system of coordinates in
 $\bR^{d}$ coincides with the one described before Lemma
\ref{lemma 2.14.1} and with the help of the mappings  
$x(y)$ and $y(x)$ introduced  in that lemma,
for $r>0$ and $z\in\bR^{d}$,
we construct  
$$
\chB^{+}_{r}=x\big(B^{+}_{r}\big), \quad 
\chB^{+}_{r}(z)=x\big(B^{+}_{r}(y(z)\big).
$$ 
By Remark \ref{remark 8,12.1} we have $\chB^{+}_{r}\subset\Omega$
for $r\leq 4\rho_{0} (\Omega)$ and
$\chB^{+}_{r}(z)\subset \chB^{+}_{4\rho_{0} (\Omega)}$
if   $ \rho>0$, 
$\rho+ r\leq 4\rho_{0}(\Omega)$,
and $z\in\chB^{+}_{\rho}$. 
Generally, these are objects in $\bR^{d}$. 

Then set
$$
\chC^{+ }_{R}=[0,R^{2})\times\chB^{+}_{R} 
$$
and for $\rho,r\leq 2\rho_{0}(\Omega)$ and $z=(t,x)$
such that $x\in \chB^{+}_{\rho}$ and $t\in\bR$ define
$$
\chC^{+}_{r}(z)=[t,t+r^{2})\times\chB^{+}_{r}(x).
$$ 
Finally, whenever it makes sense, for
  $\rho\leq\rho_{0}(\Omega)$,
and $z\in\chC^{+ }_{2\rho_{0}(\Omega)}$ introduce
\begin{equation}
                                                  \label{08,12.4}
 h^{\vsharp}_{\Omega,\gamma,\rho}(z)=
\sup\Big\{I_{r }(h,z_{0}):z_{0}\in\bR^{+}\times\Omega,
r\in(0,\rho],\chC^{+}_{r}(z_{0})\ni z \Big\},
\end{equation}
where
$$
I_{r}(h,z_{0})=\bigg(
\dashint_{ \chC^{+}_{r}(z_{0})}\dashint_{\chC^{+}_{r}(z_{0})}
|h(z_{1})-h(z_{2})|^{\gamma}\,dz_{1}dz_{2}\bigg)^{1/\gamma}.
$$
The reader should pay attention to the above curved
sharp symbol, reminding of curved boundaries.

Observe that, if $r\leq \rho\leq\rho_{0}(\Omega)$ and $z\in
\chC^{+ }_{2\rho_{0}(\Omega)}
\cap \chC^{+}_{r}(z_{0})$, then   
$$
\chC^{+}_{r}(z_{0})
\subset \chC^{+ }_{4\rho_{0}(\Omega)},
$$
 so that 
$h^{\vsharp}_{\Omega,\gamma,\rho}(z)$ is well defined on 
$\chC^{+ }_{2\rho_{0}(\Omega)}$
even if $h$ is given only on $\chC^{+ }_{4\rho_{0}(\Omega)}$
($\subset[0,16\rho_{0}^{2}(\Omega)\times\Omega$).

\begin{theorem}
                                     \label{theorem 8,12.1}
If $p\in(1,\infty)$, $\kappa\in(0,1/2]$, $0<R
\leq 2\rho_{0}(\Omega)$, 
 and $h\in L_{p}\big(\chC^{+ }_{R(1+2\kappa)}\big)$,
then
\begin{equation}
                                            \label{8,12.2}
\|h\|_{L_{p}\big( \chC^{+ }_{R}\big)}
\leq N\big\|h^{\vsharp}_{\Omega, \gamma,\kappa R}\big\|
_{L_{p}\big( \chC^{+ }_{R}\big)}+ N\kappa^{-\chi_{1}}
R^{-\chi_{2}}\big\|\,|h|^{\gamma}\big
\|^{1/\gamma}_{L_{1}\big(\chC^{+ }_{R}\big)},
\end{equation}
where  $\chi_{1}=(d+2)/\gamma$, $\chi_{2}=
(d+2)(1/\gamma-1/p )$
 and the constants
 $N$ depend only on $d,\gamma$, and $p$.

\end{theorem}

This theorem is derived from Theorem \ref{theorem 8,4.1}
by changing variables and even extension
of the functions involved across the plane $\{x^{1}=0\}$.

{\bf Acknowledgment}. Part of the work on the article was
done during the author's stay at the
Center of Mathematical Analysis and Application of
Harvard University
for one month in November 2015, and it is my great pleasure
to thank S.-T. Yau for his kind invitation.
  The author is also sincerely grateful
to the referee for very careful reading and
many comments and suggestions.

\end{document}